\documentclass{amsart}

\usepackage{amsmath}
\usepackage{amscd}
\usepackage{amssymb}

\newcommand{\bC}{{\mathbb C}}
\newcommand{\bF}{{\mathbb F}}

\newcommand{\bP}{{\mathbb P}}

\newcommand{\bZ}{{\mathbb Z}}

\newcommand{\cE}{{\mathcal E}}

\newcommand{\cP}{{\mathcal P}}

\newcommand{\cW}{{\mathcal W}}

\newtheorem{theorem}{Theorem}[section]
\newtheorem{theorem/definition}{Theorem/Definition}[section]

\newtheorem{proposition}{Proposition}[section]
\newtheorem{lemma}{Lemma}[section]

\newtheorem{corollary}{Corollary}[section]

\theoremstyle{remark}
\newtheorem{remark}{Remark}[section]

\theoremstyle{definition}
 \newtheorem{example}{Example}[section]

\begin{document}

\title
{Curve Counting and Instanton Counting}
\author{Jian Zhou}
\address{Department of Mathematical Sciences\\Tsinghua University\\Beijing, 100084, China}
\email{jzhou@math.tsinghua.edu.cn}
\begin{abstract}
We prove some combinatorial results required for
the proof of the following conjecture of Nekrasov:
The generating function of closed string invariants
in local Calabi-Yau geometries obtained by appropriate fibrations
of $A_N$ singularities over $\bP^1$
reproduce the generating function of equivariant $\hat{A}$-genera  of moduli space
of instants on $\bC^2$.
\end{abstract}
\maketitle

\section{Introduction}

Intuitively, in four-dimensional gauge theories one counts instantons,
and in two-dimensional string theory one counts holomorphic curves.
These correspond to the instanton counting and curve counting in the title.
There is no obvious connection between them,
but a remarkable physical idea called geometric engineering \cite{Kat-Kle-Vaf, Kat-May-Vaf, Kle-Ler-May-Vaf-War}
dictates that
they are indeed related.
More precisely,
the generating series of suitable integrals on the moduli spaces
of instantons and holomorphic curves respectively should be related.
By localization calculations of the equivariant $\hat{A}$-genera of the Gieseker compactifications
of moduli spaces of $SU(N)$ instantons on $\bC^2$,
Nekrasov \cite{Nek} was led to the following partition function:
\begin{eqnarray} \label{eqn:Nek}
Z = \sum_{\mu^1, \dots, \mu^N}
\varphi^{\sum_{i=1}^N |\mu^i|}
\prod_{l,n=1}^N\prod_{i,j=1}^{\infty}
\frac{\sinh \frac{\beta}{2}(a_{ln} + \hbar(\mu^l_i-\mu^n_i+j-i))}
{\sinh \frac{\beta}{2}(a_{ln} + \hbar(j-i))}.
\end{eqnarray}
(For related work, see \cite{Nak-Yos, Nek-Oko} and references therein.)
He conjectured $Z$ is the topological string partition functions
of local Calabi-Yau geometries of appropriate fibrations of the ALE spaces over $\bP^1$.

This conjecture has been studied from a physical approach in \cite{Iqb-Kas1, Iqb-Kas2}
and the more recent \cite{Hol-Iqb-Vaf}.
These works are based on the physical method of evaluating the partition functions in
local Calabi-Yau geometries by diagrammatic method developed
in \cite{Aga-Mar-Vaf, Dia-Flo-Gra, Iqb, Aga-Kle-Mar-Vaf}.
Such a method expresses the partition functions in terms of Chern-Simons link invariants,
hence leads to the rich combinatorics of representation theories of
symmetric groups and Kac-Moody algebras,
and the theory of symmetric functions.
Inspired by the results in the physical literature,
the author \cite{Zho3} has developed a mathematically rigorous approach to
the computation of partition function based on the formula for two-partition Hodge integrals,
 recently proved in joint work with C.-C. Liu and K.~Liu \cite{LLZ}.
The mathematical theory of the more general topological vertex \cite{Aga-Kle-Mar-Vaf}
is being developed jointly with J.~Li, C.-C.~Liu and K.~Liu \cite{LLLZ}.

Given the recent progresses which provides mathematical foundation for the
diagrammatic method,
it becomes possible to mathematically prove Nerasov's conjecture.
The purpose of this note is to present proofs of the combinatorial results required
for a mathematical proof of Nekrasov's conjecture
as outlined in \cite{Iqb-Kas1, Iqb-Kas2}.
Given two partitions $\mu^1$ and $\mu^2$,
define
\begin{eqnarray*}
&& K_{\mu^1\mu^2}(Q) = \sum_{\nu} Q^{|\nu|} \cW_{\mu^1\nu}\cW_{\nu\mu^2}.
\end{eqnarray*}
Iqbal and Kashani-Poor (cf. \cite[(32)]{Iqb-Kas1}) conjectured
\begin{eqnarray*}
K_{\mu^1\mu^2}(Q)
& = & K_{(0)(0)}(Q)\cW_{\mu^1}\cW_{\mu^2}
\exp \left( \sum_{n=1}^{\infty} \frac{Q^n}{n}f_{\mu^1\mu^2}(q^n)\right).
\end{eqnarray*}
for some functions $f_{\mu^1\mu^2}(q)$.
We will prove this formula in Theorem \ref{thm:K}.
Furthermore,
there is a series expansion of the form:
$$f_{\mu^1\mu^2}(q) = \sum_k C_k(\mu^1,\mu^2) q^k.$$
We will prove $C_k(\mu^1, \mu^2)$ is nonnegative.
Furthermore,
Iqbal and Kashini-Poor \cite{Iqb-Kas1} conjectured the following identity:
\begin{equation*}
\begin{split}
\prod_k (1-q^kQ)^{-2C_k(\mu^1, (\mu^2)^t)}
= & Q^{-|\mu^1|-|\mu^2|}2^{-2(|\mu^1|+|\mu^2|)}q^{-\frac{1}{2}(\kappa_{\mu^1}-\kappa_{\mu^2})}\\
& \cdot \prod_{l\neq n}\prod_{i,j\geq 1}
\frac{\sinh \frac{\beta}{2}(a_{ln}+\hbar(\mu^i_l-\mu^n_j+j-i))}
{\sinh \frac{\beta}{2}(a_{ln}+\hbar(j-i))}.
\end{split}
\end{equation*}
We will prove this identity in Corollary \ref{cor:IK2}.
Our proof is based on standard Schur function calculus.
These results and the results in \cite{Zho3} completes the proof
of Nekrasov's conjecture in the $SU(2)$ case.
For the $SU(N)$ case,
we will also prove the corresponding combinatorial assumptions made in \cite{Iqb-Kas2}.
As mentioned above the mathematical theory of the topological vertex required for
this case is being developed.
During the final stage of the preparation of this paper,
the author becomes aware of a recent  paper \cite{Egu-Kan} which has some overlaps
with this work.

The rest of this paper is arranged as follows.
We prove some results in Section \ref{sec:Prod} on infinite products associated to partitions.
We prove some summation results for skew Schur functions in \ref{sec:Schur}.
After recalling in Section \ref{sec:W} some results on $\cW_{\mu^1\mu^2}$
we study $f_{\mu^1\mu^2}(q)$ in Section \ref{sec:fmumu}.
We study $K_{\mu^1\mu^2}(Q,q)$ in Section \ref{sec:Kmumu} and Section \ref{sec:Kmumut},
and consider some generalizations in Section \ref{sec:Gen}.
In the final Section \ref{sec:Nek} we prove the $SU(2)$ case of
Nekrasov's conjecture as formulated in \cite{Iqb-Kas1}.
We also discuss the proof in the general case.

{\bf Acknowledgements}.
{\em The authors thanks Professor Kefeng Liu for suggesting the problem and
for bringing most of the references to his attentions.
He also thanks Professor Amer Iqbal and Professor Cumrun Vafa for helpful explanations.
This research is
partly supported by research grants from NSFC and Tsinghua
University.}

\section{Preliminaries Results on Partitions and infinite product formulas}
\label{sec:Prod}

\subsection{Partitions}

A partition $\mu$ is a sequence of nonincreasing nonnegative integers $\mu_1 \geq \mu_2 \geq \cdots$,
finite of which are nonzero.
The number of nonzero $\mu_i$'s in $\mu$ is called the length of $\mu$ and will be denoted by
$l(\mu)$.
The degree of $\mu$ is defined by:
$$|\mu|: = \sum_i \mu_i.$$
One can represent a partition by its Young diagram.
Given a partition $\mu$, by transposing its Young diagram one get another partition $\mu^t$.
Define
$$\kappa_{\mu} = \sum_{i=1}^{l(\mu)} \mu_i(\mu_i-2i+1).$$
There is one exceptional case:
$$\kappa_{(0)} = 0.$$
The integer $\kappa_{\nu}$ has the following nontrivial symmetric property:
$$\kappa_{\mu^t} = - \kappa_{\mu}.$$
For a proof, see e.g. \cite[Proposition 2.1]{Zho1}.
Let $x$ be the square located at the $i$-th row and the $j$-th column of the Young diagram
of $\mu$.
The content of $x$ and the hook length of $x$ is defined respectively by:
\begin{align*}
c(x) & = j -i, & h(x) & = \nu_i + \nu_j^t-i-j+1.
\end{align*}

\subsection{Infinite product formulas involving one partition}

\begin{lemma}
For a partition $\mu$,
one has the following identities:
\begin{eqnarray}
&& \sum_{1\leq i < j < \infty} (t^{\mu_i-\mu_j+j-i} - t^{j-i}) \nonumber \\
& = & \sum_{1\leq i < j \leq l(\mu)} (t^{\mu_i-\mu_j+j-i} - t^{j-i})
- \sum_{i=1}^{l(\mu)}\sum_{v=1}^{\mu_i} t^{v-i+l(\mu)} \label{eqn:Sum11}  \\
& = & - \sum_{x \in \mu} t^{h(x)}. \label{eqn:Sum12}
\end{eqnarray}
\end{lemma}

\begin{proof}
For any $n \geq l(\mu)$ one has (cf. \cite[p. 10, (2)]{Mac}):
\begin{eqnarray} \label{eqn:Summu}
&& \sum_{x \in \mu} t^{h(x)}
+ \sum_{\substack{1 \leq i < j \leq n\\1 \leq i \leq l(\mu)}} t^{\mu_i-\mu_j+j-i}
= \sum_{i=1}^{l(\mu)} \sum_{j=1}^{\mu_i-i+n} t^j.
\end{eqnarray}
In particular,
when $n = l(\mu)$,
we have:
\begin{eqnarray*}
&& \sum_{x \in \mu} t^{h(x)}
+ \sum_{1 \leq i < j \leq l(\mu)} t^{\mu_i-\mu_j+j-i} \\
& = & \sum_{i=1}^{l(\mu)} \sum_{j=1}^{\mu_i-i+l(\mu)} t^j \\
& = & \sum_{i=1}^{l(\mu)} \sum_{j=1}^{l(\mu)-i} t^j
+ \sum_{i=1}^{l(\mu)} \sum_{j=1+l(\mu)-i}^{\mu_i-i+l(\mu)} t^j \\
& = & \sum_{1 \leq i < j \leq l(\mu)} t^{j-i}
+ \sum_{i=1}^{l(\mu)}\sum_{v=1}^{\mu_i} t^{v-i+l(\mu)}.
\end{eqnarray*}
Therefore,
\begin{eqnarray} \label{eqn:Finite}
&& \sum_{1 \leq i < j \leq l(\mu)} (t^{\mu_i-\mu_j+j-i} -  t^{j-i})
- \sum_{i=1}^{l(\mu)}\sum_{v=1}^{\mu_i} t^{v-i+l(\mu)}
= - \sum_{x \in \mu} t^{h(x)}.
\end{eqnarray}

On the other hand, note
\begin{eqnarray*}
&& \sum_{i=l(\mu)+1}^n \sum_{j=i+1}^n t^{\mu_i-\mu_j+j-i}
= \sum_{i=l(\mu)+1}^n \sum_{j=i+1}^n t^{j-i} \\
& = & \sum_{i=l(\mu)+1}^{n} \sum_{j=1}^{n-i} t^j
= \sum_{i=l(\mu)+1}^{n} \sum_{j=1}^{\mu_i-i+n} t^j.
\end{eqnarray*}
Adding
$$\sum_{i=l(\mu)+1}^n \sum_{j=i+1}^n t^{\mu_i-\mu_j+j-i}
= \sum_{i=l(\mu)+1}^{n} \sum_{j=1}^{\mu_i-i+n} t^j,$$
to (\ref{eqn:Summu}) one gets:
\begin{eqnarray*}
&& \sum_{x \in \mu} t^{h(x)}
+ \sum_{1 \leq i < j \leq n} t^{\mu_i-\mu_j+j-i}
= \sum_{i=1}^{n} \sum_{j=1}^{\mu_i-i+n} t^j.
\end{eqnarray*}
Therefore
\begin{eqnarray*}
&& \sum_{1 \leq i < j \leq n} (t^{\mu_i-\mu_j+j-i} - t^{j-i})
= - \sum_{x \in \mu} t^{h(x)}
+ \sum_{i=1}^{l(\mu)} \sum_{j=1+n-i}^{\mu_i-i+n} t^j.
\end{eqnarray*}
Taking $n \to \infty$,
one then gets:
\begin{eqnarray} \label{eqn:Infinite}
&& \sum_{1 \leq i < j < \infty} (t^{\mu_i-\mu_j+j-i} - t^{j-i})
= - \sum_{x \in \mu} t^{h(x)}.
\end{eqnarray}
The proof is completed by comparing (\ref{eqn:Finite}) with (\ref{eqn:Infinite}).
\end{proof}

Here we present another proof of (\ref{eqn:Sum11}) which we will generalize in the next subsection.
\begin{eqnarray*}
&& \sum_{1\leq i < j < \infty} (t^{\mu_i-\mu_j+j-i} - t^{j-i}) \\
& = & \sum_{1\leq i < j \leq l(\mu)} (t^{\mu_i-\mu_j+j-i} - t^{j-i})
+ \sum_{1 \leq i \leq l(\mu) < j < \infty} (t^{\mu_i+j-i} - t^{j-i}) \\
& = & \sum_{1\leq i < j \leq l(\mu)} (t^{\mu_i-\mu_j+j-i} - t^{j-i})
+ \sum_{i=1}^{l(\mu)}\sum_{j=l(\mu)+1}^{\infty} t^{\mu_i+j-i}
- \sum_{i=1}^{l(\mu)}\sum_{j=l(\mu)+1}^{\infty} t^{j-i} \\
& = & \sum_{1\leq i < j \leq l(\mu)} (t^{\mu_i-\mu_j+j-i} - t^{j-i})
+ \sum_{i=1}^{l(\mu)}\sum_{v=\mu_i+1}^{\infty} t^{v-i+l(\mu)}
- \sum_{i=1}^{l(\mu)}\sum_{v=1}^{\infty} t^{v-i+l(\mu)} \\
& = & \sum_{1\leq i < j \leq l(\mu)} (t^{\mu_i-\mu_j+j-i} - t^{j-i})
- \sum_{i=1}^{l(\mu)}\sum_{v=1}^{\mu_i} t^{v-i+l(\mu)}.
\end{eqnarray*}

\begin{corollary}
For any partition $\mu$ we have:
\begin{eqnarray}
\sum_{1\leq i < j < \infty} (\mu_i-\mu_j) \nonumber
& = & \sum_{1\leq i < j \leq l(\mu)} (\mu_i-\mu_j)
- \sum_{i=1}^{l(\mu)}\sum_{v=1}^{\mu_i} (v-i+l(\mu)) \label{eqn:Sum111}  \\
& = & - \sum_{x \in \mu} h(x)
=  -\frac{1}{2}\kappa_{\mu} + |\mu|. \label{eqn:Sum121}
\end{eqnarray}
\end{corollary}

\begin{proof}
The first two identities are obtained from (\ref{eqn:Sum11}) and (\ref{eqn:Sum12}) by taking
derivative in $t$ then setting $t=1$.
To get the last equality,
note
\begin{eqnarray*}
\sum_{x \in \nu} c(x) & = &  \sum_{i=1}^{l(\nu)} \sum_{v=1}^{\nu_i} (v-i)
= \sum_{i=1}^{l(\nu)} \frac{\nu_i(\nu_i+1)}{2} -\sum_{i=1}^{l(\nu)} i\nu_i \\
& = & \frac{1}{2}\sum_{i=1}^{l(\nu)} \nu_i(\nu_i-2i+1)
= \frac{1}{2}\kappa_{\nu}.
\end{eqnarray*}
Hence
\begin{eqnarray*}
&& \sum_{1\leq i < j \leq l(\mu)} (\mu_i-\mu_j)
- \sum_{i=1}^{l(\mu)}\sum_{v=1}^{\mu_i} (v-i+l(\mu))  \\
& = & \sum_{i=1}^{l(\mu)} (l(\mu) -i)\mu_i - \sum_{j=1}^{l(\mu)} (j-1)\mu_j
- \frac{1}{2}\kappa_{\mu} - |\mu|l(\mu) \\
& = & -\frac{1}{2}\kappa_{\mu} + |\mu|.
\end{eqnarray*}
\end{proof}

As a straightforward consequence of the above cancellation arguments,
one also has:

\begin{proposition} \label{prop:Prod1}
Let $f$ be a function defined on the set of integers,
such that
$f(n) \neq 0$ for $n \neq 0$.
Then for a partition $\mu$,
one has
\begin{eqnarray}
&& \prod_{1 \leq i < j < \infty} \frac{f(\mu_i-\mu_j+j-i)}{f(j-i)} \nonumber \\
& = & \prod_{1 \leq i < j \leq l(\mu)}
\frac{f(\mu_i-\mu_j+j-i)}{f(j-i)}
\prod_{i=1}^{l(\mu)} \prod_{v=1}^{\mu_i} \frac{1}{f(v-i+l(\mu))} \\
& = & \prod_{x \in \mu} \frac{1}{f(h(x))}.
\end{eqnarray}
\end{proposition}

One can interchange $i$ with $j$ and change $t$ to $-t$ in (\ref{eqn:Sum11}) and (\ref{eqn:Sum12})
to get:
\begin{eqnarray*}
&& \sum_{1\leq j < i < \infty} (t^{\mu_i-\mu_j+j-i} - t^{j-i})  \\
& = & \sum_{1\leq j < i \leq l(\mu)} (t^{\mu_i-\mu_j+j-i} - t^{j-i})
- \sum_{i=1}^{l(\mu)}\sum_{v=1}^{\mu_i} t^{-(v-i+l(\mu))} \\
& = & - \sum_{x \in \mu} t^{-h(x)}.
\end{eqnarray*}
Hence it is easy to see:
\begin{eqnarray*}
&& \sum_{1\leq i , j < \infty} (t^{\mu_i-\mu_j+j-i} - t^{j-i})  \\
& = & \sum_{1\leq i,  j \leq l(\mu)} (t^{\mu_i-\mu_j+j-i} - t^{j-i})
- \sum_{i=1}^{l(\mu)}\sum_{v=1}^{\mu_i} t^{v-i+l(\mu)}
- \sum_{i=1}^{l(\mu)}\sum_{v=1}^{\mu_i} t^{-(v-i+l(\mu))} \\
& = & - \sum_{x \in \mu} t^{h(x)} - \sum_{x \in \mu} t^{-h(x)}.
\end{eqnarray*}
For any function $f$ satisfying the condition in Proposition \ref{prop:Prod1} one has
\begin{eqnarray*}
&& \prod_{1 \leq i , j < \infty} \frac{f(\mu_i-\mu_j+j-i)}{f(j-i)} \nonumber \\
& = & \prod_{1 \leq i, j \leq l(\mu)}
\frac{f(\mu_i-\mu_j+j-i)}{f(j-i)}
\prod_{i=1}^{l(\mu)} \prod_{v=1}^{\mu_i} \frac{1}{f(v-i+l(\mu))f(-(v-i+l(\mu)))} \\
& = & \prod_{x \in \mu} \frac{1}{f(h(x))f(-h(x))}.
\end{eqnarray*}
In particular,
when $f$ is even or odd,
\begin{eqnarray}
&& \prod_{1 \leq i , j < \infty} \frac{f(\mu_i-\mu_j+j-i)}{f(j-i)} \nonumber \\
& = & \prod_{1 \leq i, j \leq l(\mu)}
\frac{f(\mu_i-\mu_j+j-i)}{f(j-i)}
\prod_{i=1}^{l(\mu)} \prod_{v=1}^{\mu_i} \frac{1}{f(v-i+l(\mu))^2} \\
& = & \prod_{x \in \mu} \frac{1}{f(h(x))^2}.
\end{eqnarray}

\subsection{Infinite product formula involving two partitions}
In this subsection we generalize the above results.

\begin{lemma} \label{lm:Sum2}
For two partitions $\mu^1$ and $\mu^2$,
one has the following identities:
\begin{equation} \label{eqn:Summumu}
\begin{split}
& \sum_{i, j \geq 1}   (t^{\mu^1_i-\mu^2_j+j-i} - t^{j-i}) \\
= &  \sum_{i=1}^{l(\mu^1)} \sum_{j=1}^{l(\mu^2)} (t^{\mu^1_i-\mu^2_j+j-i} - t^{j-i})
-  \sum_{i=1}^{l(\mu^1)} \sum_{v=1}^{\mu^1_i} t^{v-i+l(\mu^2)}
- \sum_{j=1}^{l(\mu^2)} \sum_{v=1}^{\mu^2_j}t^{-(v-j+l(\mu^1))}.
\end{split} \end{equation}
Furthermore,
\begin{eqnarray*}
&& \sum_{i=1}^{l(\mu^1)} \sum_{j=1}^{l(\mu^2)} (\mu^1_i-\mu^2_j)
- \sum_{i=1}^{l(\mu^1)} \sum_{v=1}^{\mu^1_i} (v-i+l(\mu^2))
+ \sum_{j=1}^{l(\mu^2)} \sum_{v=1}^{\mu^2_j} (v-j+l(\mu^1)) \\
& = & -\frac{1}{2}\kappa_{\mu^1} + \frac{1}{2}\kappa_{\mu^2}.
\end{eqnarray*}
\end{lemma}

\begin{proof}
By the method in last subsection we have
\begin{eqnarray*}
&& \sum_{i, j \geq 1} (t^{\mu^1_i-\mu^2_j+j-i} - t^{j-i}) \\
& = & \sum_{i=1}^{l(\mu^1)} \sum_{j=1}^{l(\mu^2)} (t^{\mu^1_i-\mu^2_j+j-i} - t^{j-i})
+  \sum_{i=1}^{l(\mu^1)} \sum_{j = l(\mu^2)+1}^{\infty} (t^{\mu^1_i+j-i} - t^{j-i}) \\
&& +  \sum_{i = l(\mu^1)+1}^{\infty} \sum_{j=1}^{l(\mu^2)} (t^{-\mu^2_j+j-i} - t^{j-i}) \\
& = & \sum_{i=1}^{l(\mu^1)} \sum_{j=1}^{l(\mu^2)} (t^{\mu^1_i-\mu^2_j+j-i} - t^{j-i})
-  \sum_{i=1}^{l(\mu^1)} \sum_{v=1}^{\mu^1_i} t^{v-i+l(\mu^2)}
- \sum_{j=1}^{l(\mu^2)} \sum_{v=1}^{\mu^2_j}t^{-(v-j+l(\mu^1))}.
\end{eqnarray*}
Now
\begin{eqnarray*}
&&  \sum_{i=1}^{l(\nu)} \sum_{v=1}^{\nu_i} (v-i)
= \sum_{i=1}^{l(\nu)} \frac{\nu_i(\nu_i+1)}{2} -\sum_{i=1}^{l(\nu)} i\nu_i \\
& = & \frac{1}{2}\sum_{i=1}^{l(\nu)} \nu_i(\nu_i-2i+1)
= \frac{1}{2}\kappa_{\nu}.
\end{eqnarray*}
Hence
\begin{eqnarray*}
&& \sum_{i=1}^{l(\mu^1)} \sum_{j=1}^{l(\mu^2)} (\mu^1_i-\mu^2_j)
- \sum_{i=1}^{l(\mu^1)} \sum_{v=1}^{\mu^1_i} (v-i+l(\mu^2))
+ \sum_{j=1}^{l(\mu^2)} \sum_{v=1}^{\mu^2_j} (v-j+l(\mu^1)) \\
& = & |\mu^1|l(\mu^2) - |\mu^2|l(\mu^1)
- \frac{1}{2}\kappa_{\mu^1} - |\mu^1|l(\mu^2)
+ \frac{1}{2}\kappa_{\mu^2} + |\mu^2|l(\mu^1) \\
& = & -\frac{1}{2}\kappa_{\mu^1} + \frac{1}{2}\kappa_{\mu^2}.
\end{eqnarray*}
\end{proof}

\begin{proposition} \label{prop:Prod2}
For any two partitions $\mu^1$ and $\mu^2$,
and for any function $f$ satisfying the conditions in Proposition \ref{prop:Prod1},
we have
\begin{eqnarray*}
&& \prod_{i, j \geq 1} \frac{f(\mu^1_i-\mu^2_j+j-i)}{f(j-i)} \\
& = &  \prod_{i=1}^{l(\mu^1)} \prod_{j=1}^{l(\mu^2)} \frac{f(\mu^1_i-\mu^2_j+j-i)}{f(j-i)}
\prod_{i=1}^{l(\mu^1)} \prod_{v=1}^{\mu^1_i} \frac{1}{f(v-i+l(\mu^1))}
\prod_{j=1}^{l(\mu^2)} \prod_{i=1}^{\mu^2_j} \frac{1}{f(-(v-j+l(\mu^2)))}.
\end{eqnarray*}
\end{proposition}

\begin{proposition} \label{prop:Sum22}
For any two partitions $\mu^1$ and $\mu^2$,
there are integers $m_k$, $k =1 \dots, |\mu^1|+|\mu^2|$,
such that:
\begin{eqnarray} \label{eqn:Sum22}
&& \sum_{i, j \geq 1}   (t^{\mu^1_i-\mu^2_j+j-i} - t^{j-i})
= - \sum_{k=1}^{|\mu^1|+|\mu^2|} t^{m_k}.
\end{eqnarray}
\end{proposition}

\begin{proof}
We rewrite (\ref{eqn:Summumu}) as follows:
\begin{eqnarray*}
&& \sum_{i, j \geq 1}   (t^{\mu^1_i-\mu^2_j+j-i} - t^{j-i}) \\
& = &  - \sum_{i=1}^{l(\mu^1)} \sum_{j=1}^{l(\mu^2)} t^{j-i}
-  \sum_{i=1}^{l(\mu^1)} \sum_{v=1}^{\mu^1_i} t^{v-i+l(\mu^2)}
- \sum_{j=1}^{l(\mu^2)} \sum_{v=1}^{\mu^2_j}t^{-(v-j+l(\mu^1))} \\
&& + \sum_{i=1}^{l(\mu^1)} \sum_{j=1}^{l(\mu^2)} t^{\mu^1_i-\mu^2_j+j-i}.
\end{eqnarray*}
Write the three sums in the second line as $A$, $B$, and $C$, respectively.
For each pair of indices $(i, j)$ satisfying $1 \leq i \leq l(\mu^1)$
and $1 \leq j \leq l(\mu^2)$,
notice:
\begin{eqnarray*}
-(\mu^2_j + l(\mu^1)-i) < \mu_i^1 - \mu^2_j < \mu^1_i + l(\mu^2) - j.
\end{eqnarray*}
Consider three disjoint cases.
When
$$-(\mu^2_j + l(\mu^1)-i) < \mu_i^1 - \mu^2_j \leq - (1 + l(\mu^1) -i),$$
we have
$$-(\mu^2_j + l(\mu^1)-j) < \mu_i^1 - \mu^2_j + j - i \leq - (1 + l(\mu^1) -j),$$
hence $t^{\mu_i^1 - \mu^2_j + j - i}$ is cancelled by a term in $C$
of the form $t^{-(v_i-j+l(\mu^1))}$ for some $1 \leq v_i \leq \mu^2_j$.
Indeed,
$$v_i = -\mu^1_i + \mu^2_j + i - l(\mu^1).$$
Suppose for $1 \leq i_1 < i_2 \leq l(\mu^1)$,
then we have
$$v_{i_1} - v_{i_2} = -(\mu^1_{i_1} - \mu^2_{i_2} + i_1 - i_2 \leq i_1 - i_2 < 0.
$$
In other words, for different $i$, we have different $v_i$.

When
$$ 1 + l(\mu^2) - j \leq  \mu_i^1 - \mu^2_j < \mu^1_i + l(\mu^2) - j,$$
we have
$$ 1 + l(\mu^2) - i \leq  \mu_i^1 - \mu^2_j < \mu^1_i + l(\mu^2) - i,$$
then $t^{\mu_i^1 - \mu^2_j + j - i}$ is cancelled by a term in $B$
of the form $t^{\tilde{v}_j-i+l(\mu^2)}$ for some $1 \leq \tilde{v}_j \leq \mu_i^1$.
The same argument as above shows for different $j$ one has different $\tilde{v}_j$.

When
$$- (l(\mu^1) -i) \leq \mu_i^1 - \mu^2_j \leq l(\mu^2) - j,$$
we have
$$j - l(\mu^1)  \leq \mu_i^1 - \mu^2_j + j - i \leq l(\mu^2) - i.$$
If  $\mu^1_i - \mu^2_j \leq 0$,
then we say $(i, j)$ is of type $I$,
and define
\begin{align*}
\tilde{i} & =  i - (\mu^1_i - \mu^2_j), &
\tilde{j} & =  j.
\end{align*}
Otherwise we say $(i, j)$ is of type II,
and define
\begin{align*}
\tilde{i} & =  i, &
\tilde{j} & =
j + \mu^1_i - \mu^2_j.
\end{align*}
In both case we have
\begin{align*}
1 & \leq \tilde{i} \leq l(\mu^1), &
1 & \leq \tilde{j} \leq l(\mu^2),
\end{align*}
and so $t^{\mu^1_i-\mu^2_j+j-i}$ is cancelled by a term in $A$ of the form
$t^{\tilde{j}-\tilde{i}}$.
Now we show when $(i_1, j_1) \neq (i_2, j_2)$,
then $(\tilde{i}_1, \tilde{j}_1) \neq (\tilde{i}_2, \tilde{j}_2)$.
Suppose not.
When $(i_1, j_1)$ and $(i_2, j_2)$ are both of type I,
then we have
\begin{align*}
i_1 - (\mu^1_{i_1} - \mu^2_{j_1}) & = i_2 - (\mu^1_{i_2} - \mu^2_{j_2}), &
j_1 & = j_2.
\end{align*}
Suppose $i_1 < i_2$,
then we get a contradiction:
$$0 >  i_1 - i_2 = \mu^1_{i_1} - \mu^1_{i_2} \geq 0.$$
When both $(i_1, j_1)$ and $(i_2, j_2)$ are of type II,
the proof is similar.
Suppose $(i_1, j_1)$ is of type I and $(i_2, j_2)$ is of type II.
Then we have:
\begin{align*}
i_1 - (\mu^1_{i_1} - \mu^2_{j_1}) & = i_2 , &
j_1 & = j_2 +(\mu^1_{i_2} - \mu^2_{j_2}),
\end{align*}
and
\begin{align*}
\mu^1_{i_1} - \mu^2_{j_1} & \leq 0, &
\mu^1_{i_2} - \mu^2_{j_2} & > 0.
\end{align*}
Now we have $i_1 \leq i_2$, $j_1 > j_2$,
and so
\begin{align*}
\mu^1_{i_1} & \geq \mu^1_{i_2}, & \mu^2_{j_1} \leq \mu^2_{j_2}.
\end{align*}
This leads to a contradiction;
$$0 \geq \mu^1_{i_1} - \mu^2_{j_1} \geq \mu^1_{i_2} - \mu^2_{j_2}  > 0.
$$
The proof is complete.
\end{proof}

\begin{remark}
It would be interesting to find combinatorial interpretations of the integer $m_k$'s
in (\ref{eqn:Sum22}),
in the fashion of (\ref{eqn:Sum12}).
\end{remark}

\section{Preliminary Results on Skew Schur Functions}
\label{sec:Schur}

\subsection{Schur functions and skew Schur functions}

The set of Schur functions $\{s_{\mu}|\mu \in \cP\}$ form a basis
of the space of symmetric functions.
The {\em skew Schur function} $s_{\mu/\nu}$ is defined to be the symmetric function
such that
$$\langle s_{\mu/\nu}, s_{\rho}\rangle
= \langle s_{\mu}, s_{\nu}s_{\rho}\rangle.$$
Equivalently,
suppose
$$s_{\nu}s_{\rho} = \sum_\mu c^{\mu}_{\nu\rho} s_{\mu},$$
then
$$s_{\mu/\nu} = \sum_{\rho} c^{\mu}_{\nu\rho}s_{\rho}.$$
Note $s_{\mu/\nu}$ is homogeneous of degree $|\mu|-|\nu|$.
I.e.,
$$s_{\mu/\nu}(Qx) = Q^{|\mu|-|\nu|}s_{\mu/\nu}(x),$$
where $x=(x_1, x_2, \dots)$,
$Qx = (Qx_1, Qx_2, \dots)$.

\subsection{Summation formulas for Skew Schur functions}
Recall he following identity (cf. \cite[p. 93, (1)]{Mac}):
\begin{eqnarray} \label{eqn:SchurSum}
&& \sum_{\eta} s_{\eta/\mu}(x)s_{\eta/\nu}(y)
= \prod_{i, j\geq 1} (1-x_iy_j)^{-1} \cdot \sum_{\tau} s_{\mu/\tau}(y)s_{\nu/\tau}(x),
\end{eqnarray}
where $x=(x_1, x_2, \dots)$, $y=(y_1, y_2, \dots)$.
In particular,
\begin{eqnarray}
&& \sum_{\eta} s_{\eta}(x)s_{\eta}(y)
= \prod_{i, j\geq 1} (1-x_iy_j)^{-1}, \label{eqn:SchurSum1} \\
&& \sum_{\eta} s_{\eta/\mu}(x)s_{\eta}(y)
= \prod_{i, j\geq 1} (1-x_iy_j)^{-1} \cdot  s_{\mu}(y).  \label{eqn:SchurSum2}
\end{eqnarray}

\begin{lemma}
The following identity holds:
\begin{equation} \label{eqn:SchurSum3}
\begin{split}
& \sum_{\nu^1, \dots, \nu^{N}}\sum_{\eta^1, \dots, \eta^{N-1}}
\prod_{k=1}^{N} s_{\nu^k/\eta^{k-1}}(x^k)Q_k^{|\nu^k|} s_{\nu^k/\eta^{k}}(y^k) \\
= & \prod_{1 \leq k < l \leq N+1} \prod_{i, j \geq 1} (1- Q_kQ_{k+1} \cdots Q_{l-1} x^k_iy^{l-1}_j)^{-1},
\end{split} \end{equation}
where $\eta^0=\eta^N = (0)$,
$x=(x^k_1, x^k_2, \dots)$, $y^k=(y^k_1, y^k_2, \dots)$.
\end{lemma}

\begin{proof}
We repeated use (\ref{eqn:SchurSum}) and its special cases (\ref{eqn:SchurSum1}) and (\ref{eqn:SchurSum2}):
\begin{eqnarray*}
&& \sum_{\nu^1, \dots, \nu^{N}}\sum_{\eta^1, \dots, \eta^{N-1}}
\prod_{k=1}^{N} s_{\nu^k/\eta^{k-1}}(x^k)Q_k^{|\nu^k|} s_{\nu^k/\eta^{k}}(y^k) \\
& = & \sum_{\eta^1, \dots, \eta^{N-1}}
\prod_{k=1}^N Q_k^{|\eta^{k}|}
\sum_{\nu^k} s_{\nu^{k}/\eta^{k-1}}(x^k)s_{\nu^k/\eta^{k}}(Q_ky^k) \\
& = & \sum_{\eta^1, \dots, \eta^{N-1}}
\prod_{k=1}^N Q_k^{|\eta^{k}|}
\sum_{\tau^{k-1}} s_{\eta^{k}/\tau^{k-1}}(x^k)s_{\eta^{k-1}/\tau^{k-1}}(Q_ky^k)
\prod_{i, j \geq 1} (1- Q_kx^k_iy^k_j)^{-1}.
\end{eqnarray*}
Since $\eta^0=\eta^N = (0)$,
so we must have $\tau^0 = \tau^{N-1} = (0)$,
hence
\begin{eqnarray*}
&& LHS \\
& = &  \sum_{\eta^1, \dots, \eta^{N-1}}
Q_1^{|\eta^1|} s_{\eta^1}(x^1)
\prod_{k=2}^{N-1}
\sum_{\tau^{k-1}} s_{\eta^{k-1}/\tau^{k-1}}(Q_ky^k)
s_{\eta^{k}/\tau^{k-1}}(x^k)  Q_k^{|\eta^{k}|} \\
&& \cdot s_{\eta^{N-1}}(y^k) \cdot \prod_{i, j \geq 1} (1- Q_kx^k_iy^k_j)^{-1} \\
& = & \sum_{\eta^1, \dots, \eta^{N-1}}  \sum_{\tau^1, \dots, \tau^{N-2}}
\prod_{k=1}^{N-1} s_{\eta^{k}/\tau^{k-1}}(x^k)  Q_k^{|\eta^{k}|} s_{\eta^{k}/\tau^{k}}(Q_{k+1}y^{k+1}) \\
&& \cdot \prod_{k=1}^N  \prod_{i, j \geq 1} (1- Q_kx^k_iy^k_j)^{-1},
\end{eqnarray*}
where $\tau^0=\tau^{N-1} = (0)$.
To summarize,
we have obtained an identity of the form:
\begin{eqnarray*}
&& \sum_{\nu^1, \dots, \nu^{N}}\sum_{\eta^1, \dots, \eta^{N-1}}
\prod_{k=1}^{N} s_{\nu^k/\eta^{k-1}}(x^k)Q_k^{|\nu^k|} s_{\nu^k/\eta^{k}}(y^k) \\
& = & \sum_{\eta^1, \dots, \eta^{N-1}}  \sum_{\tau^1, \dots, \tau^{N-2}}
\prod_{k=1}^{N-1} s_{\eta^{k}/\tau^{k-1}}(x^k)  Q_k^{|\eta^{k}|} s_{\eta^{k}/\tau^{k}}(Q_{k+1}y^{k+1}) \\
&& \cdot \prod_{k=1}^N  \prod_{i, j \geq 1} (1- Q_kx^k_iy^k_j)^{-1},
\end{eqnarray*}
where $\eta^0=\eta^N=\tau^0=\tau^{N-1} = (0)$.
Hence (\ref{eqn:SchurSum3}) follows by induction.
\end{proof}

\section{Preliminary Results on the Topological Vertex}
\label{sec:W}

\subsection{Some previous results}
Write
$q^{\rho} = (q^{-\frac{1}{2}}, q^{-\frac{3}{2}}, \dots)$.
In \cite{Zho1} we have shown:
\begin{eqnarray} \label{eqn:smu}
&& s_{\mu}(q^{-\rho})
= (-1)^{|\mu|} q^{-\kappa_{\mu}/4}V_{\mu}(q),
\end{eqnarray}
where
\begin{eqnarray} \label{eqn:Vmu}
&& V_{\mu}(q) =
\frac{1}{\prod_{e \in \mu} (q^{h(e)/2} - q^{-h(e)/2})}.
\end{eqnarray}
We have shown in \cite{Zho2}:
\begin{eqnarray} \label{eqn:snuqmu+rho}
&& s_{\nu}(q^{\mu+\rho})
= (-1)^{|\nu|} q^{\kappa_{\nu}/2} \sum_{\eta}
\frac{s_{\mu/\eta}(q^{-\rho})}{s_{\mu}(q^{-\rho})} s_{\nu/\eta}(q^{-\rho}).
\end{eqnarray}

\subsection{Symmetries}

Note we have the obvious symmetries:
\begin{eqnarray}
&& V_{\mu^t}(q) = V_{\mu}(q), \label{eqn:Symm1} \\
&& V_{\mu}(q^{-1}) = (-1)^{|\mu|} V_{\mu}(q), \label{eqn:Symm2} \\
&& V_{\mu^t}(q^{-1}) = (-1)^{|\mu|} V_{\mu}(q). \label{eqn:Symm3}
\end{eqnarray}
We also have
\begin{align}
|\mu^t| & = |\mu|,
& \kappa_{\mu^t} & = - \kappa_{\mu}. \label{eqn:Symm4}
\end{align}
We will see later (\ref{eqn:Symm1}) - (\ref{eqn:Symm4}) are responsible for
the symmetries of $\cW_{\mu}$, $\cW_{\mu, \nu}$, and $\cW_{\mu^1, \mu^2, \mu^3}$.
For this purpose we will need the following:

\begin{proposition}
We have the following symmetries:
\begin{eqnarray}
s_{\mu^t}(q^{-\rho})
& = & q^{\kappa_{\mu}/2}s_{\mu}(q^{-\rho}), \label{eqn:smuSymm1} \\
s_{\mu}(q^{\rho})
& = & (-1)^{|\mu|} q^{\kappa_{\mu}/2}s_{\mu}(q^{-\rho}), \label{eqn:smuSymm2} \\
s_{\mu^t}(q^{\rho})
& = & (-1)^{|\mu|} s_{\mu}(q^{-\rho}), \label{eqn:smuSymm3} \\
s_{\lambda/\mu}(q^{\nu+\rho}) & = & (-1)^{|\lambda|-|\mu|}s_{\lambda^t/\mu^t}(q^{-\nu-\rho}).
\label{eqn:slambda/mu}
\end{eqnarray}
\end{proposition}

\begin{proof}
By (\ref{eqn:smu}) and (\ref{eqn:Vmu}) we have
\begin{eqnarray*}
&& s_{\mu}(q^{\frac{1}{2}}, q^{\frac{3}{2}}, \dots)
= (-1)^{|\mu|} \frac{q^{-\kappa_{\mu}/4}}{\prod_{e \in \mu} (q^{h(e)/2} - q^{-h(e)/2})}.
\end{eqnarray*}

By (\ref{eqn:smu}), (\ref{eqn:Symm1}), and (\ref{eqn:Symm4}), we have
\begin{eqnarray*}
s_{\mu^t}(q^{-\rho})
& = & (-1)^{|\mu^t|} q^{-\kappa_{\mu^t}/4}V_{\mu^t}(q)
= (-1)^{|\mu|}q^{\kappa_{\mu}/4}V_{\mu} (q)
= q^{\kappa_{\mu}/2}s_{\mu}(q^{-\rho}).
\end{eqnarray*}

By (\ref{eqn:smu}) and (\ref{eqn:Symm2}) we have
\begin{eqnarray*}
s_{\mu}(q^{\rho})
& = & (-1)^{|\mu|} q^{\kappa_{\mu}/4}V_{\mu}(q^{-1})
= q^{\kappa_{\mu}/4}V_{\mu} (q)
= (-1)^{|\mu|} q^{\kappa_{\mu}/2}s_{\mu}(q^{-\rho}).
\end{eqnarray*}

By (\ref{eqn:smu}) and (\ref{eqn:Symm3}) we have
\begin{eqnarray*}
s_{\mu^t}(q^{\rho})
& = & (-1)^{|\mu^t|} q^{\kappa_{\mu^t}/4}V_{\mu^t}(q^{-1})
= q^{-\kappa_{\mu}/4}V_{\mu} (q)
= (-1)^{|\mu|} s_{\mu}(q^{-\rho}).
\end{eqnarray*}

Finally,
\begin{eqnarray*}
s_{\mu^t/\nu^t}(q^{\rho})
& = & \sum_{\eta^t} c^{\mu^t}_{\nu^t\eta^t} s_{\eta^t}(q^{\rho}) \\
& = & \sum_{\rho} c^{\mu}_{\nu\rho}
(-1)^{|\rho|} s_{\rho}(q^{-\rho})
= (-1)^{|\mu|-|\nu|} s_{\mu/\nu}(q^{-\rho}).
\end{eqnarray*}

We begin with case of $\mu=\emptyset$.
By (\ref{eqn:snuqmu+rho}) we have
\begin{eqnarray*}
s_{\nu}(q^{\mu+\rho})
& = & (-1)^{|\nu|} q^{\kappa_{\nu}/2} \sum_{\rho}
\frac{s_{\mu/\rho}(q^{-\rho})}{s_{\mu}(q^{-\rho})} s_{\nu/\rho}(q^{-\rho}) \\
& = & (-1)^{|\nu^t|} (q^{-1})^{\kappa_{\nu^t}/2} \sum_{\rho}
\frac{(-1)^{|\mu|-|\rho|} s_{\mu^t/\rho^t}(q^{\rho})}{(-1)^{|\mu|}s_{\mu^t}(q^{\rho})}
(-1)^{|\nu|-|\rho|} s_{\nu^t/\rho^t}(q^{\rho}) \\
& = & (q^{-1})^{\kappa_{\nu^t}/2} \sum_{\rho}
\frac{ s_{\mu^t/\rho^t}(q^{\rho})}{s_{\mu^t}(q^{\rho})} s_{\nu^t/\rho^t}(q^{\rho}) \\
& = & (-1)^{|\nu|} s_{\nu^t}(q^{-\mu-\rho}).
\end{eqnarray*}
So we have
\begin{eqnarray*}
s_{\lambda/\mu}(q^{\nu+\rho})
& = & \sum_{\eta} c^{\lambda}_{\mu\eta}s_{\eta}(q^{\nu+\rho})
= \sum_{\eta} c^{\lambda^t}_{\mu^t\eta^t} (-1)^{|\eta|} s_{\eta^t}(q^{-\nu-\rho}) \\
& = & (-1)^{|\lambda|-|\mu|}s_{\lambda^t/\mu^t}(q^{-\nu-\rho}).
\end{eqnarray*}
\end{proof}

\begin{remark}
(\ref{eqn:slambda/mu}) was given in \cite[(3.10)]{ORV} without proof.
\end{remark}

\subsection{Symmetries of $\cW_{\mu}$}

Recall
\begin{eqnarray}
&& \cW_{\mu} = q^{\kappa_{\mu}/4}\prod_{1 \leq i < j \leq l(\mu)}
\frac{[\mu_i - \mu_j + j - i]}{[j-i]}
\prod_{i=1}^{l(\mu)} \prod_{v=1}^{\mu_i} \frac{1}{[v-i+l(\mu)]},
\end{eqnarray}
where as usual,
$$[m] = q^{m/2} - q^{-m/2}.$$
We have shown in \cite{Zho1}:
\begin{eqnarray}
&& \cW_{\mu}(q) = q^{\frac{\kappa_{\mu}}{4}}V_{\mu}(q),
\end{eqnarray}
and so
\begin{eqnarray}
\cW_{\mu}(q)
& = & (-1)^{|\mu|} q^{\kappa_{\mu}/2}s_{\mu}(q^{-\rho}). \label{eqn:Wmu2}
\end{eqnarray}
By (\ref{eqn:smuSymm2}) we also have
\begin{eqnarray}
\cW_{\mu}(q)
& = & s_{\mu}(q^{\rho}). \label{eqn:Wmu3}
\end{eqnarray}
It is straightforward to get the following results from (\ref{eqn:Symm1}) - (\ref{eqn:Symm3}).

\begin{proposition}
We have
\begin{eqnarray*}
&& \cW_{\mu^t}(q) = q^{-\frac{\kappa_{\mu}}{2}}\cW_{\mu}(q),  \label{eqn:Wt} \\
&& \cW_{\mu}(q^{-1}) = (-1)^{|\mu|} q^{-\frac{\kappa_{\mu}}{2}}\cW_{\mu}(q), \\
&& \cW_{\mu^t}(q^{-1}) = (-1)^{|\mu|} \cW_{\mu}(q).
\end{eqnarray*}
\end{proposition}

\subsection{Symmetries of $\cW_{\mu^1, \mu^2}$}

Recall
\begin{eqnarray}
&& \cW_{\mu, \nu} = q^{|\nu|/2} \cW_{\mu} \cdot s_{\nu}(\cE_{\mu}(q,t)),
\end{eqnarray}
where\begin{eqnarray}
&& \cE_{\mu}(q,t) = \prod_{j=1}^{l(\mu)} \frac{1+q^{\mu_j-j}t}{1+q^{-j}t}
\cdot \left(1 + \sum_{n=1}^{\infty}
\frac{t^n}{\prod_{i=1}^n (q^i-1)}\right).
\end{eqnarray}

In \cite{Zho2} we have shown:
\begin{eqnarray}
&& \cW_{\mu, \nu} = q^{|\nu|/2} \cW_{\mu} \cdot s_{\nu}(\cE_{\mu}(q,t))
= \cW_{\mu}s_{\nu}(q^{\mu+\rho}).
\end{eqnarray}
Furthermore,
by (\ref{eqn:snuqmu+rho}),
\begin{eqnarray} \label{eqn:Wmunu}
&& \cW_{\mu, \nu}(q)
= (-1)^{|\mu|+|\nu|}
q^{\frac{\kappa_{\mu}+\kappa_{\nu}}{2}}
\sum_{\eta} s_{\mu/\eta}(q^{-\rho})s_{\nu/\eta}(q^{-\rho}).
\end{eqnarray}

\begin{proposition}
We have
\begin{eqnarray}
&& \cW_{\mu^1, \mu^2}(q) = \cW_{\mu^2, \mu^1}(q), \label{eqn:Wmumu12} \\
&& \cW_{(\mu^1)^t, (\mu^2)^t}(q^{-1}) = (-1)^{|\mu^1|+|\mu^2|} \cW_{\mu^1, \mu^2}(q). \label{eqn:Wtt}
\end{eqnarray}
\end{proposition}

\begin{proof}
(\ref{eqn:Wmumu12}) follows easily from (\ref{eqn:Wmunu}).
By (\ref{eqn:slambda/mu}) and (\ref{eqn:Wmunu}), we have
\begin{eqnarray*}
&& \cW_{\mu^t, \nu^t}(q^{-1}) \\
& = & (-1)^{|\mu^t|+|\nu^t|}
q^{-\frac{\kappa_{\mu^t}-\kappa_{\nu^t}}{2}} \sum_{\eta}
s_{\mu^t/\eta^t}(q^{\rho})s_{\nu^t/\eta^t}(q^{\rho})\\
& = & q^{\frac{\kappa_{\mu}+\kappa_{\nu}}{2}} \sum_{\rho}
s_{\mu/\rho}(q^{-\rho})s_{\nu/\rho}(q^{-\rho}) \\
& = & (-1)^{|\mu^t|+|\nu^t|} \cW_{\mu, \nu}(q).
\end{eqnarray*}
This prove (\ref{eqn:Wtt}).
\end{proof}

\subsection{Reformulation of the topological vertex in terms of skew Schur functions}
The topological vertex introduced in \cite{Aga-Kle-Mar-Vaf} is defined by
\begin{eqnarray} \label{eqn:TV}
\cW_{\mu^1, \mu^2, \mu^3}
= \sum_{\rho^1, \rho^3}c_{\rho^1(\rho^3)^t}^{\mu^1(\mu^3)^t}q^{\kappa_{\mu^2}/2+\kappa_{\mu^3}/2}
\frac{\cW_{(\mu^2)^t\rho^1}\cW_{\mu^2(\rho^3)^t}}{\cW_{\mu^2\emptyset}},
\end{eqnarray}
where
$$ c_{\rho^1(\rho^3)^t}^{\mu^1(\mu^3)^t}
= \sum_{\eta} c_{\eta\rho^1}^{\mu^1}c_{\eta(\rho^3)^t}^{(\mu^3)^t}.
$$

\begin{proposition} \label{prop:TVinSchur}
We have
\begin{eqnarray*}
\cW_{\mu^1, \mu^2, \mu^3}
& = & (-1)^{|\mu^2|} q^{\kappa_{\mu^3}/2}
s_{(\mu^2)^t}(q^{-\rho}) \sum_{\eta}
s_{\mu^1/\eta}(q^{(\mu^2)^t+\rho})
s_{(\mu^3)^t/\eta}(q^{\mu^2+\rho}) \\
& = & (-1)^{|\mu^1|+|\mu^2|+|\mu^3|} q^{\kappa_{\mu^3}/2}
s_{(\mu^2)^t}(q^{-\rho}) \sum_{\eta}
s_{(\mu^1)^t/\eta}(q^{-(\mu^2)^t-\rho})
s_{\mu^3/\eta}(q^{-\mu^2-\rho}).
\end{eqnarray*}
\end{proposition}

\begin{proof}
By (\ref{eqn:TV}) and (\ref{eqn:Wmunu}) we have
\begin{eqnarray*}
&& \cW_{\mu^1, \mu^2, \mu^3} \\
& = & q^{\kappa_{\mu^2}/2+\kappa_{\mu^3}/2} \sum_{\rho^1, \rho^3,
\eta}c_{\eta\rho^1}^{\mu^1}c_{\eta(\rho^3)^t}^{(\mu^3)^t}
\frac{\cW_{(\mu^2)^t\rho^1}\cW_{\mu^2(\rho^3)^t}}{\cW_{\mu^2\emptyset}} \\
& = & q^{\kappa_{\mu^2}/2+\kappa_{\mu^3}/2}
\sum_{\eta} \frac{\sum_{\rho^1}
c_{\eta\rho^1}^{\mu^1}\cW_{(\mu^2)^t\rho^1}\cdot \sum_{\rho^3}
c_{\eta(\rho^3)^t}^{(\mu^3)^t}\cW_{\mu^2(\rho^3)^t}}{\cW_{\mu^2}} \\
& = & q^{\kappa_{\mu^2}/2+\kappa_{\mu^3}/2}
\sum_{\eta}  \frac{\sum_{\rho^1}
c_{\eta\rho^1}^{\mu^1}\cW_{(\mu^2)^t}s_{\rho^1}(q^{(\mu^2)^t+\rho})\cdot \sum_{\rho^3}
c_{\eta(\rho^3)^t}^{(\mu^3)^t}\cW_{\mu^2} s_{(\rho^3)^t}(q^{\mu^2+\rho})}{\cW_{\mu^2}} \\
& = & q^{\kappa_{\mu^2}/2+\kappa_{\mu^3}/2}
\cW_{(\mu^2)^t} \sum_{\eta}
s_{\mu^1/\eta}(q^{(\mu^2)^t+\rho})
s_{(\mu^3)^t/\eta}(q^{\mu^3+\rho}) \\
& = & q^{\kappa_{\mu^2}/2+\kappa_{\mu^3}/2} (-1)^{|(\mu^2)^t|}
q^{\kappa_{(\nu^2)^t}}s_{(\mu^2)^t}(q^{-\rho}) \sum_{\eta}
s_{\mu^1/\eta}(q^{(\mu^2)^t+\rho})
s_{(\mu^3)^t/\eta}(q^{\mu^2+\rho}) \\
& = & (-1)^{|\mu^2|} q^{\kappa_{\mu^3}/2}
s_{(\mu^2)^t}(q^{-\rho}) \sum_{\eta}
s_{\mu^1/\eta}(q^{(\mu^2)^t+\rho})
s_{(\mu^3)^t/\eta}(q^{\mu^2+\rho}) \\
& = & (-1)^{|\mu^1|+|\mu^2|+|\mu^3|} q^{\kappa_{\mu^3}/2}
s_{(\mu^2)^t}(q^{-\rho}) \sum_{\eta}
s_{(\mu^1)^t/\eta}(q^{-(\mu^2)^t-\rho})
s_{\mu^3/\eta}(q^{-\mu^2-\rho}).
\end{eqnarray*}
In the last equality we have used (\ref{eqn:slambda/mu}).
\end{proof}

\begin{remark}
Identity of this type appeared earlier
in \cite[(3.15)]{ORV} and \cite[(4.19)]{Hol-Iqb-Vaf} without proof.
\end{remark}

\subsection{Relationship between $\cW_{\mu, \nu}$ and $\cW_{\mu^1, \mu^2, \mu^3}$}

\begin{proposition} \label{prop:W(0)}
We have
\begin{eqnarray*}
\cW_{\mu^1, \mu^2, (0)}
& = & q^{\kappa_\mu^2/2} \cW_{\mu^1, (\mu^2)^t}, \\
\cW_{\mu^1, (0), \mu^3}
& = & (-1)^{|\mu^1|+|\mu^3|} q^{\kappa_{\nu^1}}\cW_{\nu^1, (\nu^3)^t}(q^{-1})
= q^{\kappa_{\nu^1}}\cW_{(\nu^1)^t, \nu^3}(q), \\
\cW_{(0), \mu^2, \mu^3}
& = & q^{\kappa_{\mu^3}/2} \cW_{\mu^2, (\mu^3)^t}(q).
\end{eqnarray*}
\end{proposition}

\begin{proof}
When $\mu^3 = (0)$,
\begin{eqnarray*}
\cW_{\mu^1, \mu^2, (0)}
& = & (-1)^{|\mu^2|} s_{(\mu^2)^t}(q^{-\rho}) s_{\mu^1}(q^{(\mu^2)^t+\rho}) \\
& = & q^{\kappa_{\mu^2}/2} \cW_{(\mu^2)^t} s_{\mu^1}(q^{(\mu^2)^t+\rho}) \\
& = & q^{\kappa_{\mu^2}/2} \cW_{(\mu^2)^t, \mu^1}
= q^{\kappa_{\mu^2}/2} \cW_{\mu^1, (\mu^2)^t}.
\end{eqnarray*}
When $\mu^1 = (0)$,
\begin{eqnarray*}
\cW_{(0), \mu^2, \mu^3}
& = & (-1)^{|\mu^2|} q^{\kappa_{\mu^3}/2}
s_{(\mu^2)^t}(q^{-\rho})
s_{(\mu^3)^t}(q^{\mu^2+\rho}) \\
& = & (-1)^{|\mu^2|} q^{\kappa_{\mu^3}/2}
s_{\mu^2}(q^{-\rho}) q^{\kappa_{\mu^2}/2}
s_{(\mu^3)^t}(q^{\mu^2+\rho}) \\
& = & q^{\kappa_{\mu^3}/2} \cW_{\mu^2} s_{(\mu^3)^t}(q^{\mu^2+\rho}) \\
& = & q^{\kappa_{\mu^3}/2} \cW_{\mu^2, (\mu^3)^t}(q).
\end{eqnarray*}
When $\mu^2 = (0)$,
\begin{eqnarray*}
&& \cW_{\mu^1, (0), \mu^3} \\
& = &  q^{\kappa_{\mu^3}/2} \sum_{\eta}
\sum_{\rho^1} c_{\eta\rho^1}^{\mu^1}\cW_{\rho^1} \cdot
\sum_{\rho^3} c_{\eta(\rho^3)^t}^{(\mu^3)^t}\cW_{(\rho^3)^t} \\
& = & q^{\kappa_{\mu^3}/2} \sum_{\eta}
\sum_{\rho^1} c_{\eta\rho^1}^{\mu^1}s_{\rho^1}(q^{-\frac{1}{2}}, q^{-\frac{3}{2}}, \dots) \cdot
\sum_{\rho^3} c_{\eta(\rho^3)^t}^{(\mu^3)^t}s_{(\rho^3)^t}(q^{-\frac{1}{2}}, q^{-\frac{3}{2}}, \dots) \\
& = & q^{\kappa_{\mu^3}/2} \sum_{\eta}
s_{\mu^1/\eta}(q^{-\frac{1}{2}}, q^{-\frac{3}{2}}, \dots) \cdot
s_{(\mu^3)^t/\eta}(q^{-\frac{1}{2}}, q^{-\frac{3}{2}}, \dots) \\
& = & (-1)^{|\mu^1|+|\mu^3|} q^{\kappa_{\nu^1}}\cW_{\nu^1, (\nu^3)^t}(q^{-1}) \\
& = & q^{\kappa_{\nu^1}}\cW_{(\nu^1)^t, \nu^3}(q)
\end{eqnarray*}
\end{proof}

\begin{remark}
The result in the above Theorem is compatible with the expected cyclic symmetry
$$\cW_{\mu^1, \mu^2, \mu^3} = \cW_{\mu^2, \mu^3, \mu^1} = \cW_{\mu^3, \mu^1, \mu^2}.$$
A mathematical proof of this symmetry is not known to the author.
\end{remark}

\subsection{Symmetries of the topological vertex}

\begin{proposition}
We have
\begin{eqnarray}
&& \cW_{(\mu^1)^t, (\mu^2)^t, (\mu^3)^t}(q)
= q^{-(\kappa_{\mu^1}+\kappa_{\mu^2}+\kappa_{\mu^3})/2}
\cW_{\mu^3, \mu^2, \mu^1}(q), \\
&& \cW_{(\mu^1)^t, \mu^2, (\mu^3)^t}(q^{-1})
= (-1)^{|\mu^1|+|\mu^2|+|\mu^3|} q^{\kappa_{\mu^2}/2} \cW_{\mu^1, \mu^2, \mu^3}(q), \\
&& \cW_{\mu^3, \mu^2, \mu^1}(q) = (-1)^{|\mu^1|+|\mu^2|+|\mu^3|}
q^{(\kappa_{\mu^1}+\kappa_{\mu^3})/2}\cW_{\mu^1, (\mu^2)^t, \mu^3}(q^{-1}).
\end{eqnarray}
\end{proposition}

\begin{proof}
By Proposition \ref{prop:TVinSchur} and (\ref{eqn:smuSymm1}) we have
\begin{eqnarray*}
&& \cW_{(\mu^1)^t, (\mu^2)^t, (\mu^3)^t}(q) \\
& = & (-1)^{|(\mu^2)^t|} q^{\kappa_{(\mu^3)^t}/2}
s_{\mu^2}(q^{-\rho}) \sum_{\eta}
s_{(\mu^1)^t/\eta}(q^{\mu^2+\rho})
s_{\mu^3/\eta}(q^{(\mu^2)^t+\rho}) \\
& = & (-1)^{|\mu^2|} q^{-\kappa_{\mu^3}/2}
\cdot q^{-\kappa_{\mu^2}/2}s_{(\mu^2)^t}(q^{-\rho}) \sum_{\eta}
s_{\mu^3/\eta}(q^{(\mu^2)^t+\rho}) s_{(\mu^1)^t/\eta}(q^{\mu^2+\rho}) \\
& = & q^{-(\kappa_{\mu^1}+\kappa_{\mu^2}+\kappa_{\mu^3})/2}
\cW_{\mu^3, \mu^2, \mu^1}(q).
\end{eqnarray*}
Similarly,
\begin{eqnarray*}
&& \cW_{(\mu^1)^t, \mu^2, (\mu^3)^t}(q^{-1}) \\
& = & (-1)^{|\mu^2|} q^{-\kappa_{(\mu^3)^t}/2}
s_{(\mu^2)^t}(q^{\rho}) \sum_{\eta}
s_{(\mu^1)^t/\eta}(-q^{(\mu^2)^t-\rho})
s_{\mu^3/\eta}(q^{-\mu^2-\rho}) \\
& = & (-1)^{|\mu^2|} q^{\kappa_{\mu^3}/2} \cdot
q^{-\kappa_{\mu^2}} s_{(\mu^2)^t}(q^{-\rho}) \\
&& \cdot \sum_{\eta}
(-1)^{|(\mu^1)^t|-|\eta|} s_{(\mu^1)^t/\eta}(q^{(\mu^2)^t+\rho}) \cdot
(-1)^{|\mu^3|-|\eta|}s_{\mu^3/\eta}(q^{\mu^2+\rho}) \\
& = & (-1)^{|\mu^1|+|\mu^2|+|\mu^3|} q^{\kappa_{\mu^2}/2} \\
&& \cdot (-1)^{|\mu^2|} q^{\kappa_{\mu^3}/2}
s_{(\mu^2)^t}(q^{-\rho}) \sum_{\eta}
s_{\mu^1/\eta}(q^{(\mu^2)^t+\rho})
s_{(\mu^3)^t/\eta}(q^{\mu^2+\rho}) \\
& = & (-1)^{|\mu^1|+|\mu^2|+|\mu^3|} q^{\kappa_{\mu^2}/2} \cW_{\mu^1, \mu^2, \mu^3}(q).
\end{eqnarray*}
The third identity can be deduced from the first two identities.
\end{proof}

\section{Results on $f_{\mu^1\mu^2}$}
\label{sec:fmumu}

\subsection{Definition of $f_{\mu^1\mu^2}$}

 For two partitions $\mu^1,
\mu^2 \in \cP$, define:
\begin{eqnarray*}
f_{\mu^1 \mu^2}(q) & = &
q \left(\sum_{i \geq 1} q^{\mu^1_i-i} \sum_{j \geq 1} q^{\mu^2_j-j}
- \sum_{i \geq 1} q^{-i} \sum_{j \geq 1} q^{-j}\right) \\
& = & \sum_{i \geq 1} q^{\mu^1_i-i+\frac{1}{2}} \sum_{j \geq 1} q^{\mu^2_j-j+\frac{1}{2}}
- \sum_{i \geq 1} q^{-i+\frac{1}{2}} \sum_{j \geq 1} q^{-j+\frac{1}{2}}.
\end{eqnarray*}
Since
\begin{eqnarray*}
&& \sum_{i \geq 1} q^{-i} = \frac{q^{-1}}{1-q^{-1}} = \frac{1}{q-1},
\end{eqnarray*}
we have
\begin{eqnarray*}
&& f_{\mu^1 \mu^2}(q) \\
& = &
q \left[\sum_{i \geq 1} (q^{\mu^1_i-i} -q^{-i}) + \frac{1}{q-1}\right]
\left[\sum_{j \geq 1} (q^{\mu^2_j-j} - q^{-j}) + \frac{1}{q-1}\right]
- \frac{q}{(1-q)^2} \\
& = & \frac{q}{(q-1)^2} \left(1+ (q-1)\sum_{i \geq 1} (q^{\mu^1_i-i} -q^{-i})\right)
\left(1+ (q-1)\sum_{j \geq 1} (q^{\mu^1_j-j} -q^{-j})\right) \\
&& - \frac{q}{(1-q)^2} \\
& = &  \frac{q}{(q-1)^2} \left(1+ (q-1)\sum_{i = 1}^{l(\mu^1)} (q^{\mu^1_i-i} -q^{-i})\right)
\left(1+ (q-1)\sum_{j = 1}^{l(\mu^2)} (q^{\mu^1_j-j} -q^{-j})\right) \\
&& - \frac{q}{(1-q)^2}.
\end{eqnarray*}
The expression in last equality appeared as (34) in \cite{Iqb-Kas1},
which also states:
\begin{eqnarray}
f_{\mu^1\mu^2}(q) & = &
\frac{\cW_{\mu^1(1)}\cW_{(1)\mu^2}- \cW_{\mu^1}\cW_{(1)}^2\cW_{\mu^2}}{\cW_{\mu^1}\cW_{\mu^2}}.
\end{eqnarray}

\subsection{The special case  $\mu^2=(0)$}
Set
$$f_{\mu}(q) : = f_{\mu(0)}(q).$$
Then we have
\begin{eqnarray}  \label{eqn:f}
f_{\mu}(q) & = & \frac{q}{q-1} \sum_{i = 1}^{l(\mu)} (q^{\mu_i-i} -q^{-i}).
\end{eqnarray}
There are two ways to rewrite this expression.
First,
\begin{eqnarray*}
f_{\mu}(q) & = &   \sum_{i = 1}^{l(\mu)}q^{-(i-1)} \frac{q^{\mu_i} - 1}{q-1}
\end{eqnarray*}
and so
\begin{eqnarray} \label{eqn:fmu}
f_{\mu}(q) & = & \sum_{i=1}^{l(\mu)} \sum_{j=1}^{\mu_i} q^{j-i}
\end{eqnarray}
(this is (36) in \cite{Iqb-Kas1}).
This expression  has the following combinatorial interpretation:
\begin{eqnarray} \label{eqn:fcx}
&& f_{\mu}(q) = \sum_{x\in \mu} q^{c(x)}.
\end{eqnarray}
Here we denote a partition by its Young diagram.
Recall the {\em content} of the square $x$ at the $i$-th row and the $j$-th column is
$$c(x) = j -i.$$

From (\ref{eqn:fcx}) it is clear that one has
\begin{eqnarray}
&& f_{\mu^t}(q) = f_{\mu}(q^{-1}), \label{eqn:fmut} \\
&& f_{\mu}(1) = |\mu|, \label{eqn:fmu1} \\
&& \left. \frac{d}{d q} f_{\mu}(q)\right|_{q=1} = \sum_{x\in \mu} c(x)
= \frac{1}{2} \kappa_{\mu}. \label{eqn:d}
\end{eqnarray}

Secondly,
\begin{eqnarray} \label{eqn:f2}
f_{\mu}(q) & = & - \sum_{j=1}^{\infty} q^j \sum_{i = 1}^{l(\mu)} (q^{\mu_i-i} -q^{-i})
= - \sum_{i = 1}^{l(\mu)} \sum_{j=1}^{\infty} (q^{\mu_i+j-i} -q^{j-i}).
\end{eqnarray}

\subsection{Series expansion}

It follows from (\ref{eqn:f}) that one has
\begin{eqnarray} \label{eqn:fmumu}
&& f_{\mu^1\mu^2}(q)
=(q-2+q^{-1})f_{\mu^1}(q)f_{\mu^2}(q)+f_{\mu^1}(q)+f_{\mu^2}(q).
\end{eqnarray}
By combining (\ref{eqn:fmut})-(\ref{eqn:d}) with (\ref{eqn:fmumu}),
one can easily prove:
\begin{eqnarray}
&& f_{(\mu^1)^t(\mu^2)^t}(q) = f_{\mu^1\mu^2}(q^{-1}),  \\
&& f_{\mu^1\mu^2}(1) = |\mu^1|+|\mu^2|,  \\
&& \left. \frac{d}{d q} f_{\mu^1\mu^2}(q)\right|_{q=1}
= \frac{1}{2} (\kappa_{\mu^1} + \kappa_{\mu^2}).
\end{eqnarray}
These are given in \cite{Iqb-Kas1} without proof.

As noted in \cite{Iqb-Kas1},
since $f_{\mu}(q)$ is a Laurent polynomial in $q$ with integral coefficients,
so is $f_{\mu^1, \mu^2}(q)$ by (\ref{eqn:fmumu}).
Hence one has an expansion of the form:
\begin{eqnarray}
&& f_{\mu^1\mu^2}(q) = \sum_k C_k(\mu^1, \mu^2)q^k,
\end{eqnarray}
where the coefficients $C_k(\mu^1,\mu^2)$ are integers
which vanish except for finitely many values of $k$.

\begin{example}
It is easy to see that
\begin{eqnarray*}
&& f_{(m)(n)} = q^{-1} + (1+ q+\cdots + q^{m-2}) + (1 + q + \cdots + q^{n-2}) + q^{m+n-1}.
\end{eqnarray*}
\end{example}

\begin{theorem}
The coefficients $C_k(\mu^1, (\mu^2)^t)$ is zero or a positive integer.
Indeed,
$$f_{\mu^1(\mu^2)^t}(q) = - \sum_{i, j \geq 1} (q^{\mu^1_i-\mu^2_j+j-i} - q^{j-i}).$$
\end{theorem}

\begin{proof}
By (\ref{eqn:f}), (\ref{eqn:fmut}) and (\ref{eqn:f2}) we have
\begin{eqnarray*}
f_{\mu^1(\mu^2)^t}(q)
& = & -(q-1)f_{\mu^1}(q) \cdot (q^{-1}-1)f_{(\mu^2)^t}(q)+f_{\mu^1}(q)+f_{(\mu^2)^t}(q) \\
& = & - (q-1)f_{\mu^1}(q) \cdot (q^{-1}-1)f_{\mu^2}(q^{-1})+f_{\mu^1}(q)+f_{\mu^2}(q^{-1}) \\
& = & - \sum_{i=1}^{l(\mu^1)} (q^{\mu^1_i-i}- q^{-i}) \cdot
\sum_{j=1}^{l(\mu^2)} (q^{-\mu^2_j+j} - q^{j}) \\
&& - \sum_{i = 1}^{l(\mu^1)} \sum_{j=1}^{\infty} (q^{\mu^1_i+j-i} -q^{j-i})
- \sum_{j = 1}^{l(\mu^2)} \sum_{i=1}^{\infty} (q^{-\mu^2_j+j-i} - q^{j-i}) \\
& = & - \sum_{i=1}^{l(\mu^1)} \sum_{j=1}^{l(\mu^2)}
(q^{\mu^1_i-\mu^2_j+j-i}- q^{\mu^1_i+j-i}- q^{-\mu^2_j+j-i} + q^{j-i}) \\
&& - \sum_{i = 1}^{l(\mu^1)} \sum_{j=1}^{\infty} (q^{\mu^1_i+j-i} -q^{j-i})
- \sum_{j = 1}^{l(\mu^2)} \sum_{i=1}^{\infty} (q^{-\mu^2_j+j-i} - q^{j-i}) \\
& = & - \sum_{i=1}^{l(\mu^1)} \sum_{j=1}^{l(\mu^2)}
(q^{\mu^1_i-\mu^2_j+j-i}-  q^{j-i}) \\
&& - \sum_{i = 1}^{l(\mu^1)} \sum_{j=l(\mu^2)+1}^{\infty} (q^{\mu^1_i+j-i} -q^{j-i})
- \sum_{j = 1}^{l(\mu^2)} \sum_{i=l(\mu^1)+1}^{\infty} (q^{-\mu^2_j+j-i} - q^{j-i}) \\
& = & - \sum_{i, j \geq 1} (q^{\mu^1_i-\mu^2_j+j-i} - q^{j-i}).
\end{eqnarray*}
Here in the last equality we have used (\ref{eqn:Summumu}).
The result is proved by invoking Proposition \ref{prop:Sum22}.
\end{proof}

\subsection{Generalization}
It is straightforward  to make the following generalizations considered in \cite{Iqb-Kas1} is as follows.
We will be brief and leave it to the interested reader to check the details.

For two partitions $\mu^1,
\mu^2 \in \cP$, define:
\begin{eqnarray*}
f_{\mu^1 \mu^2}(q_1, q_2) & = &
\sqrt{q_1q_2} \left(\sum_{i \geq 1} q_1^{\mu^1_i-i} \sum_{j \geq 1} q_2^{\mu^2_j-j}
- \sum_{i \geq 1} q_1^{-i} \sum_{j \geq 1} q_2^{-j}\right).
\end{eqnarray*}
Then we have
\begin{eqnarray*}
&& f_{\mu^1 \mu^2}(q_1, q_2) \\
& = &  \frac{\sqrt{q_1q_2}}{(q_1-1)(q_2-1)}
\left(1+ (q_1-1)\sum_{i = 1}^{l(\mu^1)} (q_1^{\mu^1_i-i} -q_1^{-i})\right)
\left(1+ (q_2-1)\sum_{j = 1}^{l(\mu^2)} (q_2^{\mu^1_j-j} -q_2^{-j})\right) \\
&& - \frac{\sqrt{q_1q_2}}{(q_1-1)(q_2-1)} \\
& = & \frac{\cW_{\mu^1(1)}(q_1)\cW_{(1)\mu^2}(q_2)
- \cW_{\mu^1}(q_1)\cW_{(1)}(q_1)\cW_{(1)}(q_2)\cW_{\mu^2}(q_2)}
{\cW_{\mu^1}(q_1)\cW_{\mu^2}(q_2)}.
\end{eqnarray*}
The expression in last equality appeared in \cite[(82)]{Iqb-Kas1}.
It follows that one has
\begin{eqnarray*}
&& f_{\mu^1\mu^2}(q_1, q_2)
=\frac{(q_1-1)(q_2-1)}{\sqrt{q_1q_2}}
f_{\mu^1}(q_1)f_{\mu^2}(q_2)
+\sqrt{\frac{q_2}{q_1}}\frac{q_1-1}{q_2-1}f_{\mu^1}(q_1)
+\sqrt{\frac{q_1}{q_2}}\frac{q_2-1}{q_1-1}f_{\mu^2}(q_2).
\end{eqnarray*}

As noted in \cite{Iqb-Kas1},
$f_{\mu^1, \mu^2}(q_1, q_2)$ has an expansion of the form:
\begin{eqnarray}
&& f_{\mu^1\mu^2}(q_1, q_2) = \sum_{k_1, k_2 \in \frac{1}{2}+\bZ}
C_{k_1, k_2}(\mu^1, \mu^2)q_1^{k_1}q_2^{k_2},
\end{eqnarray}
where the coefficients $C_{k_1, k_2}(\mu^1,\mu^2)$ are integers.

\begin{theorem}
One has
$$f_{\mu^1(\mu^2)^t}(q_1, q_2)
= -  \sqrt{\frac{q_1}{q_2}} \sum_{i=1}^{\infty} \sum_{j=1}^{\infty}
(q_1^{\mu^1_i-i}q_2^{-\mu^2_j+j} - q_1^{\mu^1_i-i}q_2^j).$$
\end{theorem}

\begin{proof}We have
\begin{eqnarray*}
&& f_{\mu^1(\mu^2)^t}(q_1, q_2) \\
& = & - \sqrt{\frac{q_2}{q_1}} (q_1-1)f_{\mu^1}(q_1)
\cdot (q_2^{-1}-1)f_{(\mu^2)^t}(q_2)+f_{\mu^1}(q)\\
&& + \sqrt{\frac{q_2}{q_1}}\frac{q_1-1}{q_2-1}f_{\mu^1}(q_1)
+ \sqrt{\frac{q_1}{q_2}}\frac{q_2-1}{q_1-1} f_{(\mu^2)^t}(q) \\
& = & -  \sqrt{\frac{q_2}{q_1}} (q_1-1)f_{\mu^1}(q_1)
\cdot (q_2^{-1}-1)f_{\mu^2}(q_2^{-1}) \\
&& + \sqrt{\frac{q_1}{q_2}}\frac{(q_1-1)}{q_1(1-q_2^{-1})} f_{\mu^1}(q_1)
+ \sqrt{\frac{q_1}{q_2}}\frac{q_2^{-1}-1}{q^{-1}_2(1- q_1)} f_{\mu^2}(q_2^{-1}) \\
& = & -  \sqrt{\frac{q_1}{q_2}}\sum_{i=1}^{l(\mu^1)} (q_1^{\mu^1_i-i}- q_1^{-i})
\cdot \sum_{j=1}^{l(\mu^2)} (q_2^{-\mu^2_j+j} - q_2^{j}) \\
&& +  \sqrt{\frac{q_1}{q_2}} \frac{1}{1-q_2^{-1}}
\sum_{i = 1}^{l(\mu^1)} (q_1^{\mu^1_i-i} -q_1^{-i})
+  \sqrt{\frac{q_1}{q_2}}\frac{1}{1-q_1}
\sum_{j = 1}^{l(\mu^2)} (q_2^{-\mu^2_j+j} - q_2^{j}) \\
& = & -  \sqrt{\frac{q_1}{q_2}} \sum_{i=1}^{l(\mu^1)} \sum_{j=1}^{l(\mu^2)}
(q_1^{\mu^1_i-i}q_2^{-\mu^2_j+j}- q_1^{\mu^1_i-i}q_2^j- q_1^{-i}q_2^{-\mu^2_j+j} + q_1^{-i}q_2^j) \\
&& -  \sqrt{\frac{q_1}{q_2}} \sum_{i = 1}^{l(\mu^1)} \sum_{j=1}^{\infty}
(q_1^{\mu^1_i-i}q_2^j -q_1^{-i}q_2^j)
-  \sqrt{\frac{q_1}{q_2}} \sum_{j = 1}^{l(\mu^2)} \sum_{i=1}^{\infty}
(q_1^{-i}q_2^{-\mu^2_j+j} - q_1^{-i}q_2^j) \\
& = & -  \sqrt{\frac{q_1}{q_2}} \left( \sum_{i=1}^{\infty} \sum_{j=1}^{\infty}
-  \sum_{i = l(\mu^1)+1}^{\infty} \sum_{j=l(\mu^2)+1}^{\infty}\right)
(q_1^{\mu^1_i-i}q_2^{-\mu^2_j+j} - q_1^{-i}q_2^j) \\
& = & -  \sqrt{\frac{q_1}{q_2}} \sum_{i=1}^{\infty} \sum_{j=1}^{\infty}
(q_1^{\mu^1_i-i}q_2^{-\mu^2_j+j} - q_1^{-i}q_2^j).
\end{eqnarray*}
\end{proof}

\section{Results on $K_{\mu^1\mu^2}(Q)$}

\label{sec:Kmumu}

Following \cite{Iqb-Kas1},
define
\begin{eqnarray*}
&& K_{\mu^1\mu^2}(Q) = \sum_{\nu} Q^{|\nu|} \cW_{\mu^1\nu}\cW_{\nu\mu^2}.
\end{eqnarray*}

The main result of this section is the following ansatz conjectured by Iqbal and Kashani-Poor
(cf. \cite[(32)]{Iqb-Kas1}):

\begin{theorem} \label{thm:K}
The following identities holds:
\begin{eqnarray}
K_{\mu^1\mu^2}(Q)
& = & K_{(0)(0)}(Q)\cW_{\mu^1}\cW_{\mu^2}
\exp \left( \sum_{n=1}^{\infty} \frac{Q^n}{n}f_{\mu^1\mu^2}(q^n)\right) \\
& = & K_{(0)(0)}(Q)\cW_{\mu^1}\cW_{\mu^2}\prod_k (1-q^kQ)^{-C_k(\mu^1, \mu^2)}.
\end{eqnarray}
\end{theorem}

\begin{proof}
We have
\begin{eqnarray*}
K_{\mu^1\mu^2}(Q)
& = & \sum_{\nu} Q^{|\nu|} \cW_{\mu^1} s_{\nu}(q^{\mu^1+\rho})
\cW_{\mu^2} s_{\nu}(q^{\mu^2+\rho}) \\
& = & \cW_{\mu^1}  \cW_{\mu^2} \sum_{\nu}
s_{\nu}(Q^{\frac{1}{2}}q^{\mu^1+\rho}) s_{\nu}(Q^{\frac{1}{2}}q^{\mu^2+\rho}) \\
& = & \cW_{\mu^1}  \cW_{\mu^2} \prod_{i, j\geq 1}
\frac{1}{(1-Qq^{\mu^1_i-i+\frac{1}{2}}q^{\mu^2_j-j+\frac{1}{2}})}.
\end{eqnarray*}
Now
\begin{eqnarray*}
&& \log \prod_{i, j\geq 1}
\frac{1}{(1-Qq^{\mu^1_i-i+\frac{1}{2}}q^{\mu^2_j-j+\frac{1}{2}})} \\
& = & - \sum_{i,j\geq 1} \log(1-Qq^{\mu^1_i-i+\frac{1}{2}}q^{\mu^2_j-j+\frac{1}{2}}) \\
& = & \sum_{i, j\geq 1} \sum_{n \geq 1}
\frac{1}{n} (Qq^{\mu^1_i-i+\frac{1}{2}}q^{\mu^2_j-j+\frac{1}{2}})^n \\
& = & \sum_{n \geq 1} \frac{1}{n}
(Qq)^n \sum_{i \geq 1} q^{n(\mu^1_i-i)} \sum_{j \geq 1} q^{n(\mu^2_j-j)} \\
\end{eqnarray*}
In particular,
when $\mu^1=\mu^2=(0)$,
\begin{eqnarray*}
&& \log \prod_{i, j\geq 1}
\frac{1}{(1-Qq^{-i+\frac{1}{2}}q^{-j+\frac{1}{2}})} \\
& = & \sum_{n \geq 1} \frac{1}{n}
(Qq)^n \sum_{i \geq 1} q^{-ni} \sum_{j \geq 1} q^{-nj} \\
& = & \sum_{n \geq 1} \frac{(Qq^{-1})^n}{n(1-q^{-n})^2},
\end{eqnarray*}
hence
\begin{eqnarray} \label{eqn:K00}
K_{(0)(0)}(Q) & = & \exp\left( \sum_{n \geq 1} \frac{(Qq^{-1})^n}{n(1-q^{-n})^2}\right).
\end{eqnarray}
Now we have
\begin{eqnarray*}
&& \log \prod_{i, j\geq 1}
\frac{(1-Qq^{-i+\frac{1}{2}}q^{-j+\frac{1}{2}})}
{(1-Qq^{\mu^1_i-i+\frac{1}{2}}q^{\mu^2_j-j+\frac{1}{2}})} \\
& = & \sum_{i,j\geq 1} \log(1-Qq^{\mu^1_i-i+\frac{1}{2}}q^{\mu^2_j-j+\frac{1}{2}})
- \sum_{i,j\geq 1} \log(1-Qq^{-i+\frac{1}{2}}q^{-j+\frac{1}{2}}) \\
& = & \sum_{n \geq 1} \frac{1}{n}
(Qq)^n \left(
\sum_{i \geq 1} q^{n(\mu^1_i-i)} \sum_{j \geq 1} q^{n(\mu^2_j-j)}
- \sum_{i \geq 1} q^{-ni} \sum_{j \geq 1} q^{-nj} \right) \\
& = & \sum_{n \geq 1} \frac{1}{n}
Q^n f_{\mu^1,\mu^2}(q^n).
\end{eqnarray*}
Hence
\begin{eqnarray*}
K_{\mu^1\mu^2}(Q)
& = & \cW_{\mu^1}  \cW_{\mu^2} K_{(0)(0)}
\prod_{i, j\geq 1}
\frac{(1-Qq^{-i+\frac{1}{2}}q^{-j+\frac{1}{2}})}
{(1-Qq^{\mu^1_i-i+\frac{1}{2}}q^{\mu^2_j-j+\frac{1}{2}})} \\
& = & \cW_{\mu^1}  \cW_{\mu^2} K_{(0)(0)}
\exp \left(\sum_{n \geq 1} \frac{Q^n}{n} f_{\mu^1\mu^2}(q^n)\right) \\
& = & \cW_{\mu^1}  \cW_{\mu^2} K_{(0)(0)}
\exp \left(\sum_{n \geq 1} \frac{Q^n}{n}
\sum_k C_k(\mu^1, \mu^2)q^{kn} \right) \\
& = & \cW_{\mu^1}  \cW_{\mu^2} K_{(0)(0)}
\prod_k (1-q^kQ)^{-C_k(\mu^1,\mu^2)}.
\end{eqnarray*}
\end{proof}

\section{Results on $K_{\mu^1(\mu^2)^t}(Q)$}
\label{sec:Kmumut}

The method in last section can be modified to prove the following result.
We will show below that it is equivalent to the ansatz (73) in \cite{Iqb-Kas1}.

\begin{theorem} \label{thm:Kt}
The following identities holds:
\begin{eqnarray} \label{eqn:Kt}
\prod_k (1-q^kQ)^{-C_k(\mu^1,(\mu^2)^t)}
& = & \frac{K_{\mu^1(\mu^2)^t}(Q)}{K_{(0)(0)}(Q)\cW_{\mu^1}(q) \cW_{(\mu^2)^t}(q)} \\
& = & Q^{-(|\mu^1|+|\mu^2|)/2} 2^{-(|\mu^1|+|\mu^2|)}
q^{-\frac{1}{4}(\kappa_{\nu^1} - \kappa_{\nu^2})} \\
&& \cdot  \prod_{i, j\geq 1}
\frac{\sinh \frac{\beta}{2}(2a+\hbar(\mu^1_i-\mu^2_j+j-i))}{\sinh \frac{\beta}{2}(2a+\hbar(j-i))}.
\nonumber
\end{eqnarray}
where
$Q=e^{-2a}$, $q = e^{-\beta \hbar}$.
\end{theorem}

\begin{proof}
By (\ref{eqn:Wtt}) and (\ref{eqn:Wt})
we have
\begin{eqnarray*}
K_{\mu^1(\mu^2)^t}(Q)
& = & \sum_{\nu} Q^{|\nu|} \cW_{\mu^1\nu}(q) (-1)^{|\mu^2|+|\nu|} \cW_{\nu^t\mu^2}(q^{-1}) \\
& = & \sum_{\nu} Q^{|\nu|} \cW_{\mu^1}(q) s_{\nu}(q^{\mu^1+\rho})
(-1)^{|\mu^2|+|\nu|} \cW_{\mu^2}(q^{-1}) s_{\nu^t}(q^{-(\mu^2+\rho)}) \\
& = & \cW_{\mu^1}(q)  \cW_{\mu^2}(q) \sum_{\nu}
s_{\nu}(Q^{\frac{1}{2}}q^{\mu^1+\rho}) s_{\nu^t}(-Q^{\frac{1}{2}}q^{-(\mu^2+\rho)}) \\
& = & \cW_{\mu^1}(q)  \cW_{(\mu^2)^t}(q) \prod_{i, j\geq 1}
(1-Qq^{\mu^1_i-i+\frac{1}{2}}q^{-(\mu^2_j-j+\frac{1}{2})}) \\
& = & \cW_{\mu^1}(q) \cW_{(\mu^2)^t}(q) \prod_{i, j\geq 1} (1-Qq^{\mu^1_i-\mu^2_j+j-i}).
\end{eqnarray*}
In particular,
\begin{eqnarray*}
K_{(0)(0)}(Q) & = & \prod_{i, j\geq 1}  (1-Qq^{j-i}),
\end{eqnarray*}
hence
\begin{eqnarray*}
\frac{K_{\mu^1(\mu^2)^t}(Q)}{K_{(0)(0)}(Q)}
& =& \cW_{\mu^1}(q)  \cW_{(\mu^2)^t}(q) \prod_{i, j\geq 1}
\frac{(1-Qq^{\mu^1_i-\mu^2_j+j-i})}{(1-Qq^{j-i})}.
\end{eqnarray*}
Now we have
\begin{eqnarray*}
&&  \prod_{i, j\geq 1}
\frac{(1-Qq^{\mu^1_i-\mu^2_j+j-i})}{(1-Qq^{\mu^1_i-\mu^2_j+j-i})}
= \prod_{i, j\geq 1}\frac{(1-e^{-2a\beta-\beta\hbar(\mu^1_i-\mu^2_j+j-i)})}{(1-e^{-2a\beta-\beta\hbar(j-i)})} \\
& = & \prod_{i=1}^{l(\mu^1)} \prod_{j=1}^{l(\mu^2)}
\frac{1 - e^{-\beta(2a+\hbar(\mu^1_i-\mu^2_j+j-i))}}{1 - e^{-\beta(2a+\hbar(j-i))}}
\prod_{i=1}^{l(\mu^1)} \prod_{v=1}^{\mu^1_i} \frac{1}{1-e^{-\beta(2a+\hbar(v-i+l(\mu^1)))}} \\
&& \cdot \prod_{j=1}^{l(\mu^2)} \prod_{i=1}^{\mu^2_j} \frac{1}{1-e^{-\beta(2a-\hbar(v-j+l(\mu^2)))}} \\
& = & \prod_{i=1}^{l(\mu^1)} \prod_{j=1}^{l(\mu^2)}
\frac{\sinh\frac{\beta}{2}(2a+\hbar(\mu^1_i-\mu^2_j+j-i))}{\sinh \frac{\beta}{2}(2a+\hbar(j-i))} \\
&& \cdot \prod_{i=1}^{l(\mu^1)} \prod_{v=1}^{\mu^1_i} \frac{1}{\sinh\frac{\beta}{2}(2a+\hbar(v-i+l(\mu^1)))} \\
&& \cdot \prod_{j=1}^{l(\mu^2)} \prod_{i=1}^{\mu^2_j}
\frac{1}{\sinh\frac{\beta}{2}(2a-\hbar(v-j+l(\mu^2)))}
\cdot (e^{\beta a}/2)^{(|\mu^1|+|\mu^2|)}  \\
&& \cdot e^{-\frac{1}{2}\beta\hbar\left[\sum_{i=1}^{l(\mu^1)} \sum_{j=1}^{l(\mu^2)} (\mu^1_i-\mu^2_j)
- \sum_{i=1}^{l(\mu^1)} \sum_{v=1}^{\mu^1_i} (v-i+l(\mu^1))
+ \sum_{j=1}^{l(\mu^2)} \sum_{i=1}^{\mu^2_j} (v-j+l(\mu^2))\right]} \\
& = & Q^{-(|\mu^1|+|\mu^2|)/2} 2^{-(|\mu^1|+|\mu^2|)}
q^{-\frac{1}{4}(\kappa_{\nu^1} - \kappa_{\nu^2})}
 \prod_{i, j\geq 1}
\frac{\sinh \frac{\beta}{2}(2a+\hbar(\mu^1_i-\mu^2_j+j-i))}{\sinh \frac{\beta}{2}(2a+\hbar(j-i))}.
\end{eqnarray*}
Here we have used Lemma \ref{lm:Sum2} and Proposition \ref{prop:Prod2}.
\end{proof}

As an easy consequence,
we have the following result conjectured in \cite{Iqb-Kas1}.

\begin{corollary} \label{cor:IK2}
For any two partitions $\mu^1$ and $\mu^2$,
one has
\begin{equation} \label{eqn:IK2}
\begin{split}
\prod_k (1-q^kQ)^{-2C_k(\mu^1, (\mu^2)^t)}
= & Q^{-|\mu^1|-|\mu^2|}2^{-2(|\mu^1|+|\mu^2|)}q^{-\frac{1}{2}(\kappa_{\mu^1}-\kappa_{\mu^2})}\\
& \cdot \prod_{l\neq n}\prod_{i,j\geq 1}
\frac{\sinh \frac{\beta}{2}(a_{ln}+\hbar(\mu^i_l-\mu^n_j+j-i))}
{\sinh \frac{\beta}{2}(a_{ln}+\hbar(j-i))},
\end{split}
\end{equation}
where $(l, n) = (1, 2)$ or $(2, 1)$, $a_{12} = - a_{21} = 2a$.
\end{corollary}

The above identity (\ref{eqn:IK2}) is the technical ansatz (73) in \cite{Iqb-Kas1}.
By applying Proposition \ref{prop:Prod1} to $f(x) = q^{\frac{x}{2}} - q^{-\frac{x}{2}}$,
one gets:
\begin{eqnarray*}
&& \prod_{1 \leq i < j < \infty} \frac{[\mu_i - \mu_j+j-i]}{[j-i]}
= \prod_{1 \leq i < j \leq l(\mu)} \frac{[\mu_i-\mu_j+j-i]}{[j-i]}
\prod_{i=1}^{l(\mu)} \prod_{v = 1}^{\mu_i} \frac{1}{[v-i+l(\mu)]},
\end{eqnarray*}
and
\begin{eqnarray*}
&& \cW_{\mu}(q) = 2^{-|\mu|}q^{\kappa_{\mu}/4}\prod_{1 \leq i < j < \infty}
\frac{\sinh \frac{\beta\hbar}{2}(\mu_i-\mu_j+j-i)}{\sinh \frac{\beta\hbar}{2}(j-i)}.
\end{eqnarray*}
These are formulas (71) and (72) in \cite{Iqb-Kas1}, respectively.
In the same fashion one gets
\begin{eqnarray*}
&& \cW_{\mu^t}(q)
= (-1)^{|\mu|} \cW_{\mu}(q^{-1}) \\
& = & (-1)^{|\mu|} 2^{-|\mu|}q^{-\kappa_{\mu}/4}\prod_{1 \leq i < j < \infty}
\frac{-\sinh \frac{\beta\hbar}{2}(\mu_i-\mu_j+j-i)}{-\sinh \frac{\beta\hbar}{2}(j-i)} \\
& = & 2^{-|\mu|}q^{-\kappa_{\mu}/4}\prod_{1 \leq i < j < \infty}
\frac{\sinh \frac{\beta\hbar}{2}(\mu_i-\mu_j+j-i)}{\sinh \frac{\beta\hbar}{2}(j-i)}
\end{eqnarray*}
Combined with (\ref{eqn:Kt}) one gets:

\begin{theorem} \label{thm:Ksinh}
The following identity holds:
\begin{equation} \label{eqn:K2}
\begin{split}
& \frac{K_{\mu^1(\mu^2)^t}(Q)^2}{K_{(0)(0)}(Q)^2} \\
= & Q^{-(|\mu^1|+|\mu^2|)} 2^{-4(|\mu^1|+|\mu^2|)}
\prod_{l, n =1}^2 \prod_{i, j\geq 1}
\frac{\sinh \frac{\beta}{2}(a_{ln}+\hbar(\mu^1_i-\mu^2_j+j-i))}{\sinh \frac{\beta}{2}(a_{ln}+\hbar(j-i))},
\end{split} \end{equation}
where $a_{11}=a_{22} = 0$,
$a_{12} = - a_{21} = 2a$.
\end{theorem}

\section{Generalizations}
\label{sec:Gen}

The results in this section are inspired by \cite[\S 4.3]{Iqb-Kas2}.
Define:
\begin{eqnarray*}
&& \tilde{K}_{\mu^1\cdots \mu^N}(Q_{1}, \dots, Q_{N-1})
= \sum_{\nu^1, \dots, \nu^{N-1}} \prod_{k=1}^{N}
q^{\kappa_{\nu^k}/2} \cW_{\nu^{k-1}\mu^k(\nu^{k})^t}(q) Q_{k}^{|\nu^k|},
\end{eqnarray*}
where $\nu^0=\nu^N = (0)$, $Q_N = 1$.
When $N=2$,
we have by Proposition \ref{prop:W(0)}
\begin{eqnarray*}
\tilde{K}_{\mu^1\cdots \mu^N}(Q)
& = & \sum_{\nu}
q^{\kappa_{\nu}/2} \cW_{(0)\mu^1\nu^t}(q) Q^{|\nu|}\cW_{\nu\mu^2(0)} \\
& = & q^{\kappa_{\mu^2}/2}\sum_{\nu} \cW_{\mu^1\nu}(q) Q^{|\nu|}\cW_{\nu(\mu^2)^t} \\
& = & q^{\kappa_{\mu^2}/2} K_{\mu^1(\mu^2)^t}.
\end{eqnarray*}

Theorem \ref{thm:K} can be generalized as follows.

\begin{theorem} \label{thm:KGen}
We have the following identity:
\begin{eqnarray*}
&& \frac{\tilde{K}_{\mu^1\cdots \mu^N}(Q_{1}, \dots, Q_{N-1})}
{\prod_{1 \leq k < l \leq N} K_{(0)(0)}(Q_k\cdots Q_{l-1})} \\
& = & \prod_{k=1}^{N} \cW_{\mu^k}(q) \cdot
\exp \left( \sum_{1 \leq k < l \leq N}
\sum_{n \geq 1} \frac{(Q_k\cdots Q_{l-1})^n}{n} f_{\mu^i(\mu^j)^t}(q^n) \right) \\
& = & \prod_{k=1}^{N} \cW_{\mu^k}(q)
\prod_{1 \leq k < l \leq N} \prod_{n \geq 1}(1 - q^nQ_k \cdots Q_{l-1})^{-C_n(\mu^k, (\mu^l)^t)}.
\end{eqnarray*}
\end{theorem}

\begin{proof}
By Proposition \ref{prop:TVinSchur}, (\ref{eqn:smuSymm3}), and (\ref{eqn:Wmu3}),
we have
\begin{eqnarray*}
\cW_{\mu^1, \mu^2, \mu^3}(q)
& = & q^{\kappa_{\mu^3}/2} \cW_{\mu^2}(q) \sum_{\eta}
s_{\mu^1/\eta}(q^{(\mu^2)^t+\rho})
s_{(\mu^3)^t/\eta}(q^{\mu^2+\rho}).
\end{eqnarray*}
and so
\begin{eqnarray*}
&& \tilde{K}_{\mu^1\cdots \mu^N}(Q_{1}, \dots, Q_{N-1}) \\
& = & \sum_{\nu^1, \dots, \nu^{N-1}} \prod_{k=1}^{N}
q^{\kappa_{\nu^k}/2} \cW_{\nu^{i-1}\mu^k(\nu^{k})^t}(q) Q_{k}^{|\nu^k|}
\hspace{.2in} (\nu^0=\nu^N = (0)) \\
& = & \sum_{\nu^1, \dots, \nu^{N-1}}
\prod_{k=1}^{N} \cW_{\mu^k}(q)
\sum_{\eta^{k-1}} s_{\nu^{k-1}/\eta^{k-1}}(q^{(\mu^k)^t+\rho})
s_{\nu^k/\eta^{k-1}}(q^{\mu^k+\rho}) Q_k^{|\nu^k|} \\
&& \hspace{2in} (\nu^0=\nu^N = \eta^0=\eta^{N-1} = (0)) \\
& = & \prod_{k=1}^{N} \cW_{\mu^k}(q)
\sum_{\nu^1, \dots, \nu^{N-1}}\sum_{\eta^1, \dots, \eta^{N-2}} \prod_{k=1}^{N-1}
s_{\nu^k/\eta^{k-1}}(q^{\mu^k+\rho}) Q_k^{|\nu^k|}
s_{\nu^k/\eta^k}(q^{(\mu^{k+1})^t+\rho}) \\
&& \hspace{2in} (\eta^0=\eta^{N-1} = (0)) \\
& = &  \prod_{k=1}^{N} \cW_{\mu^k}(q) \cdot
\prod_{1 \leq k < l \leq N} \prod_{i, j \geq 1}
(1- Q_kQ_{k+1} \cdots Q_{l-1} q^{\mu^k_i-i+\frac{1}{2}}q^{(\mu^l)^t_j-j+\frac{1}{2}})^{-1}.
\end{eqnarray*}
In the last equality we have used (\ref{eqn:SchurSum3}).
Hence
\begin{eqnarray*}
&& \frac{\tilde{K}_{\mu^1\cdots \mu^N}(Q_{1}, \dots, Q_{N-1})}
{\prod_{1 \leq k < l \leq N} K_{(0)(0)}(Q_k\cdots Q_{l-1})} \\
& = &  \prod_{k=1}^{N} \cW_{\mu^k}(q) \cdot
\prod_{1 \leq k < l \leq N} \prod_{i, j \geq 1}
\frac{(1- Q_kQ_{k+1} \cdots Q_{l-1} q^{-i+\frac{1}{2}}q^{-j+\frac{1}{2}})}
{(1- Q_kQ_{k+1} \cdots Q_{l-1} q^{\mu^k_i-i+\frac{1}{2}}q^{(\mu^l)^t_j-j+\frac{1}{2}})}  \\
& = & \prod_{k=1}^{N} \cW_{\mu^k}(q) \cdot
\exp \left( \sum_{1 \leq k < l \leq N}
\sum_{n \geq 1} \frac{(Q_k\cdots Q_{l-1})^n}{n} f_{\mu^i(\mu^j)^t}(q^n) \right) \\
& = & \prod_{k=1}^{N} \cW_{\mu^k}(q) \cdot
\prod_{1 \leq k < l \leq N} \prod_{n \geq 1}(1 - q^nQ_k \cdots Q_{l-1})^{-C_n(\mu^k, (\mu^l)^t)}.
\end{eqnarray*}
In the last two equalities we have used Theorem \ref{thm:K} and its proof.
\end{proof}

Theorem \ref{thm:Ksinh} can be generalized as follows.
Let
\begin{align*}
Q_i & = e^{-\beta a_i}, & a_{ij}&=a_i-a_j.
\end{align*}

\begin{theorem}
The following identity holds:
\begin{eqnarray*}
&& \frac{\tilde{K}_{\mu^1\cdots \mu^N}(Q_{1}, \dots, Q_{N-1})^2}
{\prod_{1 \leq k < l \leq N} K_{(0)(0)}(Q_k\cdots Q_{l-1})^2} \\
& = & 2^{-2N\sum_{i=1}^N |\mu^i|}
\prod_{k=1}^{N-1}
Q_k^{-[(N-k)(|\mu^1|+ \cdots +|\mu^k|) + k(|\mu^{k+1}| + \cdots + |\mu^N|)]} \\
&& q^{-\frac{1}{2}\sum_{k=1}^N (N-2k)\kappa_{\nu^k}}
\prod_{1 \leq k, l \leq N} \prod_{i, j\geq 1}
\frac{\sinh \frac{\beta}{2}(2a_{kl}+\hbar(\mu^k_i-\mu^l_j+j-i))}
{\sinh \frac{\beta}{2}(2a_{kl}+\hbar(j-i))}.
\end{eqnarray*}
\end{theorem}

\begin{proof}
By Theorem \ref{thm:KGen} and Theorem \ref{thm:Kt} we have
\begin{eqnarray*}
&& \frac{\tilde{K}_{\mu^1\cdots \mu^N}(Q_{1}, \dots, Q_{N-1})}
{\prod_{1 \leq k < l \leq N} K_{(0)(0)}(Q_k\cdots Q_{l-1})} \\
& = & \prod_{k=1}^{N} \cW_{\mu^k}(q) \cdot
\prod_{1 \leq k < l \leq N} \prod_{n \geq 1}(1 - q^nQ_k \cdots Q_{l-1})^{-C_n(\mu^k, (\mu^l)^t)} \\
& = & \prod_{k=1}^{N} 2^{-|\mu^k|}q^{\kappa_{\mu^k}/4}\prod_{1 \leq i < j < \infty}
\frac{\sinh \frac{\beta\hbar}{2}(\mu^k_i-\mu^k_j+j-i)}{\sinh \frac{\beta\hbar}{2}(j-i)} \\
&& \cdot \prod_{1 \leq k < l \leq N}
(Q_k \cdots Q_{l-1})^{-(|\mu^k|+|\mu^l|)/2} 2^{-(|\mu^k|+|\mu^l|)}
q^{-\frac{1}{4}(\kappa_{\nu^k} - \kappa_{\nu^l})} \\
&& \cdot  \prod_{i, j\geq 1}
\frac{\sinh \frac{\beta}{2}(2a_{kl}+\hbar(\mu^k_i-\mu^l_j+j-i))}
{\sinh \frac{\beta}{2}(2a_{kl}+\hbar(j-i))} \\
& = & 2^{-N\sum_{i=1}^N |\mu^i|}
\prod_{k=1}^{N-1}
Q_k^{-\frac{1}{2}[(N-k)(|\mu^1|+ \cdots +|\mu^k|) + k(|\mu^{k+1}| + \cdots + |\mu^N|)]} \\
&& q^{-\frac{1}{4}\sum_{k=1}^N (N-2k)\kappa_{\nu^k}}
\prod_{k=1}^{N} \prod_{1 \leq i < j < \infty}
\frac{\sinh \frac{\beta\hbar}{2}(\mu^k_i-\mu^k_j+j-i)}{\sinh \frac{\beta\hbar}{2}(j-i)} \\
&& \cdot \prod_{1 \leq k < l \leq N} \prod_{i, j\geq 1}
\frac{\sinh \frac{\beta}{2}(2a_{kl}+\hbar(\mu^k_i-\mu^l_j+j-i))}
{\sinh \frac{\beta}{2}(2a_{kl}+\hbar(j-i))}.
\end{eqnarray*}
The proof is completed by taking squares.
In the above we have implicitly used the following calculations.
\begin{eqnarray*}
&& \sum_{k=1}^N |\mu^k| + \sum_{1 \leq k < l \leq N} (|\mu^k|+|\mu^l|) \\
& = & \sum_{k=1}^N |\mu^k| + \sum_{k=1}^{N-1} \sum_{l=k+1}^N |\mu^k|
+  \sum_{k=1}^{N-1} \sum_{l=k+1}^N |\mu^l| \\
& = & \sum_{k=1}^N |\mu^k| + \sum_{k=1}^{N} (N-k) |\mu^k|
+  \sum_{l=2}^N \sum_{k=1}^{l-1} |\mu^l| \\
& = & \sum_{k=1}^N |\mu^k| + \sum_{k=1}^{N} (N-k) |\mu^k|
+  \sum_{l=1}^N (l-1) |\mu^l| \\
& = & N \sum_{k=1}^N |\mu^k|;
\end{eqnarray*}
\begin{eqnarray*}
&&  \prod_{1 \leq k < l \leq N} (Q_k \cdots Q_{l-1})^{-(|\mu^k|+|\mu^l|)/2} \\
& = & \prod_{m=1}^N Q_m^{-\frac{1}{2} \sum_{k \leq m < l \leq N} (|\mu^k|+|\mu^l|)} \\
& = &  \prod_{m=1}^N Q_m^{-\frac{1}{2} \sum_{k=1}^m \sum_{l =m+1}^N (|\mu^k|+|\mu^l|)} \\
& = & \prod_{m=1}^N
Q_m^{-\frac{1}{2} \sum_{k=1}^m (N-m) |\mu^k|+ \sum_{l=m+1}^N m|\mu^l|};
\end{eqnarray*}
\begin{eqnarray*}
&& \sum_{k=1}^{N} \kappa_{\mu^k} - \sum_{1 \leq k < l \leq N}
(\kappa_{\nu^k} - \kappa_{\nu^l}) \\
& = & \sum_{k=1}^{N} \kappa_{\mu^k} - \sum_{k=1}^{N-1} \sum_{l=k+1}^N
\kappa_{\nu^k} -  \sum_{l=2}^{N} \sum_{k=1}^{l-1} \kappa_{\nu^l} \\
& = & \sum_{k=1}^{N} \kappa_{\mu^k} - \sum_{k=1}^{N} (N-k)\kappa_{\nu^k}
+ \sum_{l=1}^{N} (l-1) \kappa_{\nu^l} \\
& = & \sum_{k=1}^N (N-2k)\kappa_{\nu^k}.
\end{eqnarray*}
\end{proof}

\section{Nekrasov Conjecture}
\label{sec:Nek}

\subsection{The $SU(2)$ case}
The partition function of local Calabi-Yau geometry of the canonical line bundle on
$\bP^1\times \bP^1$ has been calculated in
\cite{Aga-Mar-Vaf, Iqb, Aga-Kle-Mar-Vaf} by physical method,
and in \cite{Zho3} by mathematical method.
The result can be written as follows:
\begin{eqnarray*}
&& Z(Q_B, Q_F; q) \\
& = & \sum_{\mu^1, \nu^1, \mu^2, \nu^2} Q_B^{|\mu^1|+|\mu^2|}Q_F^{|\nu^1|+|\nu^2|}
\cW_{\mu^1, \nu^1}\cW_{\nu^1\mu^2}\cW_{\mu^2\nu^2}\cW_{\nu^2\mu^1} \\
& = & \sum_{\mu^1, \mu^2} Q_B^{|\mu^1|+|\mu^2|} K_{\mu^1, \mu^2}(Q_F)^2
=  \sum_{\mu^1, \mu^2} Q_B^{|\mu^1|+|\mu^2|} K_{\mu^1, (\mu^2)^t}(Q_F)^2.
\end{eqnarray*}
In general,
for a Hirzebruch sufrace $\bF_m$ ($m=0, 1,2$),
the partition function on its canonical line bundle is
\begin{eqnarray*}
Z^{(m)}(Q_B, Q_D;q) = \sum_{\mu^1, \mu^2} Q_B^{|\mu^1|+|\mu^2|}
[(-1)^{|\mu^1|+|\mu^2|}Q_F^{|\mu^2|}q^{-\frac{1}{2}(\kappa_{\mu^1}+\kappa_{\mu^2})}]^m
K_{\mu^1, (\mu^2)^t}(Q_F)^2.
\end{eqnarray*}
Now we can prove the following form of Nekrasov's conjecture formulated in \cite{Iqb-Kas1}.

\begin{theorem}
We have
\begin{equation}
\begin{split}
\frac{Z^{(m)}(Q_B, Q_F; q)}{K^2_{(0)(0)}(Q_F)}
= & \sum_{\mu^1, \mu^2} ((-1)^m Q_B2^{-4}Q_F^{-1})^{|\mu^1|+|\mu^2|}
Q_F^{m|\mu^2|}q^{-\frac{m}{2}(\kappa_{\mu^1}+\kappa_{\mu^2})} \\
& \cdot \prod_{l,n=1}^2\prod_{i, j\geq 1}
\frac{\sinh \frac{\beta}{2}(a_{ln}+\hbar(\mu^1_i-\mu^2_j+j-i))}{\sinh \frac{\beta}{2}(a_{ln}+\hbar(j-i))}.
\end{split} \end{equation}
In particular,
\begin{eqnarray*}
&& \frac{Z^{(0)}(Q_B, Q_F; q)}{K^2_{(0)(0)}(Q_F)}
= \sum_{\mu^1, \mu^2} \left(\frac{Q_B}{2^{4}Q_F}\right)^{|\mu^1|+|\mu^2|}
\prod_{l,n=1}^2\prod_{i, j\geq 1}
\frac{\sinh \frac{\beta}{2}(a_{ln}+\hbar(\mu^1_i-\mu^2_j+j-i))}{\sinh \frac{\beta}{2}(a_{ln}+\hbar(j-i))}.
\end{eqnarray*}
\end{theorem}

\begin{proof}
By the results in last section we have
\begin{eqnarray*}
&& \frac{Z^{(m)}(Q_B, Q_F; q)}{K^2_{(0)(0)}(Q_F)} \\
& = &   \sum_{\mu^1, \mu^2} Q_B^{|\mu^1|+|\mu^2|}
[(-1)^{|\mu^1|+|\mu^2|}Q_F^{|\mu^2|}q^{-\frac{1}{2}(\kappa_{\mu^1}+\kappa_{\mu^2})}]^m\\
&& \cdot
Q_F^{-(|\mu^1|+|\mu^2|)} 2^{-4(|\mu^1|+|\mu^2|)} \\
&& \prod_{l,n=1}^2\prod_{i, j\geq 1}
\frac{\sinh \frac{\beta}{2}(a_{ln}+\hbar(\mu^1_i-\mu^2_j+j-i))}{\sinh \frac{\beta}{2}(a_{ln}+\hbar(j-i))} \\
& = & \sum_{\mu^1, \mu^2} ((-1)^m Q_B2^{-4}Q_F^{-1})^{|\mu^1|+|\mu^2|}
Q_F^{m|\mu^2|}q^{-\frac{m}{2}(\kappa_{\mu^1}+\kappa_{\mu^2})} \\
&& \prod_{l,n=1}^2\prod_{i, j\geq 1}
\frac{\sinh \frac{\beta}{2}(a_{ln}+\hbar(\mu^1_i-\mu^2_j+j-i))}{\sinh \frac{\beta}{2}(a_{ln}+\hbar(j-i))}.
\end{eqnarray*}
\end{proof}

\subsection{The general case}

Suggested by (66) and (67) in \cite{Iqb-Kas2},
the partition function should be
\begin{eqnarray*}
\tilde{Z}^{(m)} & = & \sum_{\mu^1, \dots, \mu^N} (\prod_{i=1}^NQ_{b_i}^{|\mu^i|})
M^{(N-2)}(q; \mu^1, \dots, \mu^i)M^{(m)}(q; \mu^1, \dots, \mu^N) \\
&& \cdot \tilde{K}^2_{\mu^1\dots, \mu^N}(Q_1, \dots, Q_{N-1}) \\
\end{eqnarray*}
where $M^{(m)}$ is a framing factor given by:
\begin{eqnarray*}
&& M^{(m)}(q; \mu^1, \dots,\mu^N)
= (-1)^{(N+m)\sum_{i=1}^N|\mu^i|}
q^{\frac{1}{2}\sum_{i=1}^N (N+m-2i)\kappa_{\mu^i}}.
\end{eqnarray*}
For $m=0$,
one then gets:
\begin{eqnarray*}
&& \frac{\tilde{Z}^{(0)}}{\prod_{1 \leq k < l \leq N} K_{(0)(0)}(Q_k\cdots Q_{l-1})} \\
& = & \sum_{\mu^1, \dots, \mu^N} \varphi^{\sum_{i=1}^N |\mu^i|}
\prod_{1 \leq k, l \leq N} \prod_{i, j \geq 1}
\frac{\sinh \frac{\beta}{2}(2a_{kl}+\hbar(\mu^k_i-\mu^l_j+j-i))}
{\sinh \frac{\beta}{2}(2a_{kl}+\hbar(j-i))},
\end{eqnarray*}
for suitably defined $\varphi$ similar to the $N=2$ case.


\begin{thebibliography}{99}


\bibitem{Aga-Kle-Mar-Vaf}
M.~Aganagic, A.~Klemm, M.~Mari\~no, C.~Vafa,
{\em The topological vertex},
preprint, hep-th/0305132.

\bibitem{Aga-Mar-Vaf}
M.~Aganagic, M.~Mari\~no, C.~Vafa,
{\em All loop topogoical string amplitudes from Chern-Simons theory},
preprint, hep-th/0206164.



\bibitem{Dia-Flo-Gra}
D.E.~Diaconescu, B.~Florea, A.~Grassi,
{\em Geometric transitions, del Pezzosurfaces and open string instantons},
Adv. Theor. Math. Phys. {\bf 6} (2003), 643, hep-th/0206163.

\bibitem{Egu-Kan}
T. Eguchi, H. Kanno,
{\em Topological Strings and Nekrasov¡¯s formulas},
preprint, hep-th/0310235.

\bibitem{Hol-Iqb-Vaf}
T.~Holowood, A.~Iqbal, C.~Vafa,
{\em Matrix models, geometric engineering and elliptic genera},
preprint, hep-th/0310272.

\bibitem{Iqb}
A.~Iqbal,
{\em All genus topological string amplitudes and 5-brane webs as Feynman diagrams},
preprint, hep-th/0207114.

\bibitem{Iqb-Kas1}
A.~Iqbal, A.-K.~Kashani-Poor,
{\em $SU(N)$ geometries and topological string amplitudes}, preprint, hep-th/0306032.

\bibitem{Iqb-Kas2}
A.~Iqbal, A.-K.~Kashani-Poor,
{\em Instanton counting and Chern-Simons theory},
hep-th/0212279.



\bibitem{Kat-Kle-Vaf}
S. Katz, A. Klemm and C. Vafa,
{\em Geometric engineering of quantum field theories},
Nucl. Phys. B{\bf 497} (1997), 173-195, hep-th/9609239.

\bibitem{Kat-May-Vaf}
S. Katz, P. Mayr and C. Vafa,
{\em Mirror symmetry and exact solution of 4D N = 2 gauge
theories. I}, Adv. Theor. Math. Phys. {\bf 1} (1998), 53-??,  hep-th/9706110.

\bibitem{Kle-Ler-May-Vaf-War}
A.~Klemm, W.~Lerche, P.~Mayr, C.~Vafa, N.P.~Warner,
{\em Self-dual strings and $N=2$ supersymmetric field theory},
Nucl. Phys. B {\bf 477} (1996), 746-766, hep-th/9604034.



\bibitem{LLLZ}
J.~Li, C.-C.~Liu, K.~Liu, J.~Zhou,
in preparation.

\bibitem{LLZ}
C.-C.~Liu, K.~Liu, J.~Zhou,
{\em A formula on two-partition Hodge integrals},
preprint, math.AG/0310272.



\bibitem{Mac}
I.G.~MacDonald,
{\em Symmetric functions and Hall polynomials}, 2nd edition.Claredon Press, 1995.

\bibitem{Nak-Yos}
H.~Nakajima, K.~Yoshioka,
{\em Instanton counting on blowup, I}, preprint,  math.AG/0306198.


\bibitem{Nek}
N.~A.~Nekrasov,
{\em Seiberg-Witten prepotential from instanton counting},
preprint, hep-th/0206161.

\bibitem{Nek-Oko}
N.~Nekrasov, A.~Okounkov,
{\em Seiberg-Witten theory and random partitions},
preprint, hep-th/0306238.

\bibitem{ORV}
A.~Okounkov, N.~Reshtihkin, C.~Vafa,
{\em Quantum Calabi-Yau and Classical Crystals}, hep-th/0309208.

\bibitem{Zho1}
J.~Zhou,
{\em Hodge integrals, Hurwitz numbers, and symmetric groups},
preprint, math.AG/0308024.

\bibitem{Zho2}
J.~Zhou,
{\em A conjecture on Hodge integrals},
preprint, math.AG/0310282.

\bibitem{Zho3}
J.~Zhou,
{\em Localizations on moduli Spaces and free field realizations of Feynman rules},
preprint, math.AG/0310283.


\end{thebibliography}
\end{document}